\documentclass[12pt,reqno]{amsart}
\usepackage{amssymb}
\usepackage{amsxtra}
\usepackage[mathscr]{eucal}

\newcommand{\dst}{\displaystyle}
\newcommand{\e}{\varepsilon}
\newcommand{\G}{\Gamma}
\newcommand{\Z}{{\mathbf Z}}

\newcommand{\Q}{{\mathbf Q}}
\newcommand{\C}{{\mathbf C}}
\newcommand{\eop}{\hfill$\square$}
\newcommand{\us}{\{1\}}
\newcommand\shuff{
  \setlength{\unitlength}{.4pt}
  \begin{picture}(40,20)
    \put(10,2){\line(1,0){20}} \put(10,2){\line(0,1){10}}
    \put(20,2){\line(0,1){10}} \put(30,2){\line(0,1){10}}
  \end{picture}}

\theoremstyle{plain}
\newtheorem{Thm}{Theorem}

\theoremstyle{definition}

\theoremstyle{remark}

\newtheorem*{Notn}{Notation and Terminology}

\numberwithin{equation}{section}

\begin{document}

\title[Goldbach Variations]{Thirty-Two Goldbach Variations}

\date{\today}

\author[J.~M.~Borwein]{Jonathan~M. Borwein}
\address{Faculty of Computer Science\\
         Dalhousie University\\
         Halifax, Nova Scotia B3H 1W5\\
         Canada}
\email[Jonathan M.~Borwein]{jborwein@cs.dal.ca}
\urladdr{http://www.cs.dal.ca/$\sim$jborwein}
\thanks{Research of the first author supported by NSERC}

\author[D.~M.~Bradley]{David~M. Bradley}
\address{Department of Mathematics \& Statistics\\
         University of Maine\\
         5752 Neville Hall
         Orono, Maine 04469-5752\\
         U.S.A.}
\email[David M.~Bradley]{bradley@math.umaine.edu,
dbradley@member.ams.org}
\urladdr{http://www.umemat.maine.edu/faculty/bradley/}



\begin{abstract}
We give thirty-two diverse proofs of a small mathematical gem---the
fundamental Euler sum identity
$$\zeta(2,1)=\zeta(3)=8\,\zeta(\overline2,1).$$ We also discuss various
generalizations for multiple harmonic (Euler) sums and some of their
many connections, thereby illustrating  both the wide variety of
techniques fruitfully used to study such sums and the attraction of
their study.
\end{abstract}

\maketitle

\tableofcontents \interdisplaylinepenalty=500

\section{Introduction}\label{sect:Intro}

There are several ways to introduce and make attractive  a new or
unfamiliar subject. We choose to do so by emulating Glen Gould's
passion for Bach's \emph{Goldberg variations}. We shall illustrate
most of the techniques used to study Euler sums by focusing almost
entirely on the identities of~\eqref{z21} and \eqref{2bar1}
\[
   \sum_{n=1}^\infty \frac{1}{n^2}\sum_{m=1}^{n-1}\frac1{m} =
    \sum_{n=1}^\infty \frac{1}{n^3} = 8\sum_{n=1}^\infty
    \frac{(-1)^n}{n^2}\sum_{m=1}^{n-1}\frac1{m}
\]
and some of their many generalizations.

\subsection{Euler, Goldbach and the birth of
${\boldsymbol\zeta}$.} What follows is a  transcription of
correspondence between Euler and Goldbach \cite{eg1742} that led to
the origin of the zeta-function and multi-zeta values, see also
\cite{exp1, exp2, dunham}.

\begin{quote}
\noindent
 {\bf 59. Goldbach an Euler, Moskau, 24. Dez. 1742.}\footnote{AAL: F.136, Op. 2, Nr.8, Blatt 54--55.}
[\dots]\emph{Als ich neulich die vermeinten summas der beiden
letzteren serierum in meinem vorigen Schreiben wieder betrachtet,
habe ich alsofort wahrgenommen, da\ss~ selbige aus einem blo\ss em
Schreibfehler entstanden, von welchem es aber in der Tat hei\ss et:
Si non errasset, fecerat ille minus}.\footnote{Frei zitiert nach
Marcus Valerius Martialis, I, 21,9.}
\end{quote}

This is the letter in which Goldbach precisely formulates the series
which sparked Euler's further investigations into what would become
  the zeta-function.  These investigations were apparently
due to a serendipitous mistake. The above translates as follows:
\begin{quote}
\emph{When I recently considered further the indicated sums of the
last two series in my previous letter, I realized immediately that
the same series arose due to a mere writing error, from which indeed
the saying goes, ``Had one not erred, one would have achieved
less."\footnote{\rm {\sl Opera Omnia}, vol. IVA4, Birkh\"auser
Verlag.}}
\end{quote}
Goldbach continues
\begin{quote}
\emph{Ich halte daf\"ur, da\ss~ es ein problema problematum ist, die
summam huius:}
\begin{eqnarray*}
1+ \frac{1}{2^n}\left(1+ \frac{1}{2^m}\right) +
\frac{1}{3^n}\left(1+ \frac{1}{2^m} +\frac{1}{3^m}\right)
+\frac{1}{4^n}\left(1+ \frac{1}{2^m}
+\frac{1}{3^m}+\frac{1}{4^m}\right)+ etc.
\end{eqnarray*}
\emph{in den casibus zu finden, wo $m$ et $n$ nicht numeri integri
pares et sibi aequales sind, doch gibt es casus, da die summa
angegeben werden kann, exempli gr[atia], si $m=1$, $n=3$, denn es
ist}
\begin{eqnarray*}
 1+ \frac{1}{2^3}\left(1+ \frac{1}{2}\right) +  \frac{1}{3^3}\left(1+
\frac{1}{2} +\frac{1}{3}\right)  +  \frac{1}{4^3}\left(1+
\frac{1}{2} +\frac{1}{3}+\frac{1}{4}\right)+ etc. =
\frac{\pi^4}{72}.
\end{eqnarray*}
\end{quote}

\subsection{The Modern Language of Euler Sums}
For positive integers $s_1,\dots,s_m$ and signs $\sigma_j = \pm 1$,
consider~\cite{BBB} the $m$-fold Euler sum
\[
   \zeta(s_1,\dots,s_m;\sigma_1,\dots,\sigma_m)
   := \sum_{k_1>\cdots>k_m>0}\;\prod_{j=1}^m
   \frac{\sigma_j^{k_j}}{k_j^{s_j}}.
\]
As is now customary, we combine strings of exponents and signs by
replacing $s_j$ by $\overline s_j$ in the argument list if and only
if $\sigma_j=-1$, and denote $n$ repetitions of a substring $S$ by
$\{S\}^n$.  Thus, for example, $\zeta(\overline1)=-\log 2$,
$\zeta(\{2\}^3)=\zeta(2,2,2)=\pi^6/7!$ and
\begin{equation}
   \zeta(s_1,\dots,s_m) = \sum_{k_1>\cdots>k_m>0}\;\prod_{j=1}^m
   k_j^{-s_j}.
\label{mzvdef}
\end{equation}
The identity
\begin{equation}\label{z21}
   \zeta(2,1) = \zeta(3)
\end{equation}
goes back to Euler~\cite{LE1}~\cite[p.\ 228]{LE2} and has since been
repeatedly rediscovered, (see,  e.g.,
~\cite{Briggs,Bruckman,Farnum,Klamkin}). In this language Goldbach
had found $$\zeta(3,1)+\zeta(4) =\frac{\pi^4}{72}.$$

 The more general formula
\begin{equation}\label{EulerReduction}
    2\zeta(m,1)= m\zeta(m+1)-\sum_{j=1}^{m-2}\zeta(j+1)\zeta(m-j),
    \qquad 2\le m\in\Z
\end{equation}
is also due to Euler~\cite{LE1}~\cite[p.\ 266]{LE2}.
Nielsen~\cite[p.\ 229]{Niels1}~\cite[p.\ 198]{Niels2}~\cite[pp.\
47--49]{Niels3} developed a method for
obtaining~\eqref{EulerReduction} and related results based on
partial fractions.  Formula~\eqref{EulerReduction} has also been
rediscovered many
times~\cite{Williams,State,SitSarm,GP,Bracken,Vowe}.  Crandall and
Buhler~\cite{CranBuh} deduced~\eqref{EulerReduction} from their
general infinite series formula which expresses $\zeta(s,t)$  for
real $s>1$ and $t\ge 1$ in terms of Riemann zeta values.

Study of the the multiple zeta function~\eqref{mzvdef} led to the
discovery of a new generalization of~\eqref{z21}, involving nested
sums of arbitrary depth:
\begin{equation}\label{z21^n}
   \zeta(\{2,1\}^n ) = \zeta(\{3\}^n), \qquad n\in\Z^{+}.
\end{equation}
Although numerous proofs of~\eqref{z21} and~\eqref{EulerReduction}
are known (we give many in the sequel), the only proof
of~\eqref{z21^n} of which we are aware involves making a simple
change of variable in a multiple iterated integral (see \cite{BBB,
BBBLa,BowBradSurvey} and~\eqref{duality} below).

An alternating version of~\eqref{z21} is
\begin{equation}\label{2bar1}
   8\zeta(\overline{2}, 1) = \zeta(3),
\end{equation}
which has also resurfaced from time to time~\cite[p.\ 50]{Niels3}
~\cite[(2.12)]{Sub85}~\cite[p.\ 267]{Butzer} and hints at the
generalization
\begin{equation}\label{2bar1^n}
   8^n\zeta(\{\overline 2, 1\}^n) \stackrel{?}{=} \zeta(\{3\}^n),\qquad
   n\in\Z^{+},
\end{equation}
originally conjectured in~\cite{BBB}, and which still remains
open---despite abundant, even overwhelming, evidence \cite{DJD}.


 \subsection{Hilbert and Hardy Inequalities}\label{sect:hardy}
  Much of the early 20th century history---and philosophy---of the
\emph{``~`bright' and amusing"}  subject of
 inequalities charmingly discussed in G.H. Hardy's retirement lecture as
London Mathematical Society Secretary, \cite{ghh}. He
comments~\cite[p.\ 474]{ghh} that \emph{Harald Bohr is reported to
have remarked ``Most analysts spend half their time hunting
through the literature for inequalities they want to use, but
cannot prove."}

 Central to Hardy's essay are:

 \begin{Thm} {\rm (\bf Hilbert)}  For  non-negative sequences $(a_n)$
and $(b_n)$, not both zero, and for\\ $1 \le p, q \le \infty$ with
$1/p+1/q=1$ one has
\begin{eqnarray} \label{hilbert-p}\sum_{n=1}^\infty\sum_{m=1}^\infty \frac{a_n\,b_m}{n+m}
< \pi\,{\rm csc}\left(\frac{\pi}p\right) \|a_n\|_p\,\|b_n\|_q.
\end{eqnarray}
\end{Thm}
 \begin{Thm} {\rm (\bf Hardy)} For a non-negative sequence $(a_n)$ and for $p>1$
\begin{eqnarray}\label{h-ineq}\sum_{n=1}^\infty\left(\frac{a_1+a_2+\cdots
+a_n}n\right)^p \le \left(\frac{p}{p-1}\right)^p \,\sum_{n=1}^\infty
a_n^p.\end{eqnarray}
\end{Thm}
We return to these inequalities in Section \ref{sect:witten}.

Hardy~\cite[p.\ 485]{ghh} remarks that his {\emph{``own theorem
was discovered as a by-product of my own attempt to find a really
simple and elementary proof of Hilbert's."}}  He reproduces
Elliott's proof of (\ref{h-ineq}), writing ``\emph{it can hardly
be possible to find a proof more concise or elegant}" and also
``I\emph{ have given nine [proofs] in a lecture in Oxford, and
more have been found since then.}" (See~\cite[p.\ 488]{ghh}.)

\subsection{Our Motivation and Intentions}

We wish to emulate Hardy and to present proofs that are either
elementary, bright and amusing, concise or elegant--- ideally all at
the same time! In doing so we note that:

\begin{enumerate} \item $\zeta(3)$, while provably irrational,  is still quite
mysterious, see \cite{bt,agm} and \cite{exp2}. Hence, exposing more
relationships and approaches can only help. We certainly hope one of
them will lead to a proof of conjecture (\ref{2bar1^n}).
\item Identities for $\zeta(3)$ are abundant and diverse. We give
three each of which is the entry-point to a fascinating set:
\begin{itemize}
\item Our first  favourite is a \emph{binomial sum }\cite{ag} that played a role in Ap\'ery's 1976 proof,
 see \cite{agm,bt} and \cite[Chapter 3]{exp2}, of
the irrationality of $\zeta(3)$:
\begin{eqnarray}\label{z3}
\zeta(3) &=& \frac5 2\,\sum_{k=1}^{\infty} \frac{(-1)^{k+1}} {k^3\,
{2k \choose k}}.
\end{eqnarray}
\item Our second is Broadhurst's binary \emph{BBP formula} \cite{broad}:
$$\zeta(3)=\frac{48}7\,\mathcal{S}_{1}(1,-7,-1,10,-1,-7,1,0 )+\frac{32}7\,\mathcal{S}_{3}(1,1,-1,-2,-1,1,1,0),$$
where $\mathcal{S}_{p}(a_1,a_2,\ldots, a_8):=\sum_{k=1}^\infty
{a_k}2^{-\lfloor{p(k+1)/2\rfloor}}k^{-3},$ and the coefficients
$a_k$ repeat modulo 8. We refer to \cite[Chapter 3]{exp1} for the
digit properties of such formulae. Explicitly,
\begin{eqnarray*}
\zeta(3) &=& \frac{1}{672} \sum_{k=0}^\infty \frac{1}{2^{12 k}}
\left[\frac{2048}{(24 k + 1)^3} - \frac{11264}{(24 k + 2)^3} -
\frac{1024}{(24 k + 3)^3} + \frac{11776}{(24 k + 4)^3} \right. \\
&& \left. - \frac{512}{(24 k + 5)^3} + \frac{4096}{(24 k + 6)^3} +
\frac{256}{(24 k + 7)^3} + \frac{3456}{(24 k + 8)^3} +
\frac{128}{(24 k + 9)^3} \right. \\
&& \left. - \frac{704}{(24 k + 10)^3} - \frac{64}{(24 k + 11)^3} -
\frac{128}{(24 k + 12)^3} - \frac{32}{(24 k + 13)^3} -
\frac{176}{(24 k + 14)^3} \right. \\
&& \left. + \frac{16}{(24 k + 15)^3}  + \frac{216}{(24 k + 16)^3} +
\frac{8}{(24 k + 17)^3} + \frac{64}{(24 k + 18)^3} - \frac{4}{(24 k
+ 19)^3} \right. \\
&& \left. + \frac{46}{(24 k + 20)^3} - \frac{2}{(24 k + 21)^3}  -
\frac{11}{(24 k + 22)^3} + \frac{1}{(24 k + 23)^3} \right].
\end{eqnarray*}
It was this discovery that lead Bailey and Crandall to their
striking recent work on normality of BBP constants \cite[Chapter
4]{exp1}.

\item Our third favourite due to Ramanujan \cite[p.\ 138]{exp2} is
the \emph{hyperbolic series} approximation
$$\zeta  \left( 3 \right) ={\frac {7\,{\pi }^{3}}{180}}-2\,\sum _{k=1}^{\infty }{\frac {1}{{k}^{3}
 \left( {e^{2\,\pi \,k}}-1 \right) }},
$$
in which the `error' is $\zeta(3)-7\,{\pi }^{3}/180 \approx
-0.003742745$, and which to our knowledge is the `closest' one gets
to writing $\zeta(3)$ as a rational multiple of $\pi^3$.
\end{itemize}
 \item Often results about $\zeta(3)$
are more precisely results about $\zeta(2,1)$ or
$\zeta(\overline{2},1)$, as we shall exhibit.
\item Double and multiple sums are still under-studied and
under-appreciated. We should like to partially redress that.
\item One can now
prove these seemingly analytic facts in an entirely finitary manner
via words over alphabets, dispensing with notions of infinity and
convergence;
 \item  Many subjects are
touched upon---from computer algebra, integer relation methods,
generating functions and techniques of integration to
polylogarithms, hypergeometric and special functions,
non-commutative rings, combinatorial algebras and Stirling
numbers---so that most readers will find a proof worth showing in an
undergraduate class;
\item For example, there has been an explosive recent interest in q-analogues, see
\S\ref{sect:zud}, and in quantum field theory, algebraic K-theory
and knot theory, see \cite{exp1,zagier}.
\end{enumerate}

 For some of the broader issues relating
to Euler sums, we refer the reader to  the survey
articles~\cite{exp1,BowBradSurvey,Cartier,Wald,Wald2,zagier,Zud}.
Computational issues are discussed in~\cite{exp2, Cran} and to an
extent in~\cite{BBBLa}.

\subsection{Further Notation}

\begin{Notn} For positive integer $N$, denote the $N$th partial sum
of the harmonic series by $H_N := \sum_{n=1}^N 1/n$.  We also use
$\psi = \Gamma'/\Gamma$ to denote the logarithmic derivative of
the Euler gamma-function (also referred to as the \emph{digamma}
function), and recall the identity $\psi(N+1)+\gamma=H_N$, where
$\gamma=0.5772156649\ldots$ is \emph{Euler's constant.}  Where
convenient, we employ the Pochhammer symbol $(a)_n=
a(a+1)\cdots(a+n-1)$ for complex $a$ and non-negative integer $n$.
As usual, the \emph{Kronecker} $\delta_{m,n}$ is 1 if $m=n$ and
$0$ otherwise.
\end{Notn}

We  organize our proofs by technique, although clearly this is
somewhat arbitrary as many proofs fit well within more than one
category. Broadly their sophistication increases as we move through
the paper. In some of the later sections the proofs become more
schematic. We invite readers to send additional  selections for our
collection, a collection which for us has all the beauty of Blake's
grain of sand\footnote{ William Blake from  \emph{Auguries of
Innocence}.}:
\begin{center}
``\emph{To see a world in a grain of sand\\ And a heaven in a wild
flower,\\
Hold infinity in the palm of your hand\\ And eternity in an hour.}"
\end{center}

\section{Telescoping and Partial Fractions}\label{sect:telescope}
For a quick proof of~\eqref{z21}, consider
\[
   S:=  \sum_{n,k>0} \frac{1}{nk(n+k)}
   = \sum_{n,k>0}\frac1{n^2} \left(\frac1k-\frac1{n+k}\right)
   = \sum_{n=1}^\infty \frac1{n^2} \sum_{k=1}^n \frac1k
   = \zeta(3)+\zeta(2,1).
\]
On the other hand,
\[
   S = \sum_{n,k>0} \left(\frac1n+\frac1k\right)\frac{1}{(n+k)^2}
   = \sum_{n,k>0}\frac1{n(n+k)^2}+\sum_{n,k>0}\frac{1}{k(n+k)^2}
   = 2\zeta(2,1),
\]
by symmetry. \qed

The above argument goes back at least to Steinberg~\cite{Klamkin}.
See also~\cite{Klamkin3}.

For~\eqref{2bar1}, first consider
\begin{align}
   \zeta(\overline 2,\overline 1)+\zeta(3)
   &= \sum_{n=1}^\infty \frac{(-1)^n}{n^2}\sum_{k=1}^n
   \frac{(-1)^k}{k}
   = \sum_{n=1}^\infty \frac{(-1)^n}{n^2}\sum_{k=1}^\infty
     \bigg(\frac{(-1)^k}{k}-\frac{(-1)^{n+k}}{n+k}\bigg)\nonumber\\
   &=\sum_{n=1}^\infty\frac{(-1)^n}{n^2}\sum_{k=1}^\infty
     (-1)^k\bigg(\frac{n+k-(-1)^nk}{k(n+k)}\bigg)\nonumber\\
   &=\sum_{n,k>0}\frac{(-1)^{n+k}}{nk(n+k)}+\sum_{n,k>0}\frac{(-1)^{n+k}}{n^2(n+k)}
    -\sum_{n,k>0}\frac{(-1)^k}{n^2(n+k)}\nonumber\\
   &=\sum_{n,k>0}\bigg(\frac1n+\frac1k\bigg)\frac{(-1)^{n+k}}{(n+k)^2}
     +\zeta(\overline
     1,2)-\sum_{n,k>0}\frac{(-1)^n(-1)^{n+k}}{n^2(n+k)}\nonumber\\
   &=\sum_{n,k>0}\frac{(-1)^{n+k}}{n(n+k)^2}+\sum_{n,k>0}\frac{(-1)^{n+k}}{k(n+k)^2}
     +\zeta(\overline1,2)-\zeta(\overline1,\overline2)\nonumber\\
   &=2\zeta(\overline2,1)+\zeta(\overline1,2)-\zeta(\overline1,\overline2).
\label{alt3}
\end{align}

Similarly,
\begin{align}
   \zeta(2,\overline1)+\zeta(\overline3)
   &=\sum_{n=1}^\infty\frac{1}{n^2}\sum_{k=1}^n\frac{(-1)^k}{k}
    =\sum_{n=1}^\infty\frac1{n^2}\sum_{k=1}^\infty\bigg(\frac{(-1)^k}{k}-\frac{(-1)^{n+k}}{n+k}\bigg)\nonumber\\
   &=\sum_{n=1}^\infty\frac1{n^2}\sum_{k=1}^\infty(-1)^k\bigg(\frac{n+k-(-1)^nk}{k(n+k)}\bigg)\nonumber\\
   &=\sum_{n,k>0}\frac{(-1)^k}{nk(n+k)}+\sum_{n,k>0}\frac{(-1)^k}{n^2(n+k)}-\sum_{n,k>0}\frac{(-1)^{n+k}}{n^2(n+k)}\nonumber\\
   &=\sum_{n,k>0}\bigg(\frac1n+\frac1k\bigg)\frac{(-1)^k}{(n+k)^2}+\sum_{n,k>0}\frac{(-1)^n(-1)^{n+k}}{n^2(n+k)}
    -\zeta(\overline1,2)\nonumber\\
   &=\sum_{n,k>0}\frac{(-1)^n(-1)^{n+k}}{n(n+k)^2}+\sum_{n,k>0}\frac{(-1)^k}{k(n+k)^2}+\zeta(\overline1,\overline2)
    -\zeta(\overline1,2)\nonumber\\
   &=\zeta(\overline2,\overline1)+\zeta(2,\overline1)+\zeta(\overline1,\overline2)-\zeta(\overline1,2).
\label{alt3bar}
\end{align}

Adding equations~\eqref{alt3} and~\eqref{alt3bar} now gives
\begin{equation}
\label{dejavu}
   2\zeta(\overline2,1)=\zeta(3)+\zeta(\overline3),
\end{equation}
i.e.\
\[
   8\zeta(\overline2,1)=4\sum_{n=1}^\infty\frac{1+(-1)^n}{n^3}=4\sum_{m=1}^\infty
   \frac{2}{(2m)^3} = \zeta(3),
\]
which is~\eqref{2bar1}. \qed

\section{Finite Series Transformations}\label{sect:finite}
For any positive integer $N$, we have
\begin{equation}\label{z21finite}
   \sum_{n=1}^N \frac1{n^3}
  -\sum_{n=1}^N \frac1{n^2}\sum_{k=1}^{n-1}\frac1k
   =\sum_{n=1}^N \frac1{n^2}\sum_{k=1}^n\frac{1}{N-k+1}
\end{equation}
by induction.  Alternatively, consider
\[
   T := \sum_{\substack{n,k=1\\k\ne n}}^N \frac{1}{nk(k-n)}
   = \sum_{\substack{n,k=1\\k\ne n}}^N \bigg(\frac1n-\frac1k\bigg)
     \frac{1}{(k-n)^2}
   =0.
\]
On the other hand,
\begin{align*}
   T &= \sum_{\substack{n,k=1\\k\ne n}}^N \frac{1}{n^2}
        \bigg(\frac{1}{k-n}-\frac1k\bigg)\\
     &= \sum_{n=1}^N\frac{1}{n^2}\bigg(\sum_{k=1}^{n-1}\frac{1}{k-n}
      +\sum_{k=n+1}^N\frac1{k-n}-\sum_{k=1}^N
      \frac1k+\frac1n\bigg)\\
     &= \sum_{n=1}^N\frac1{n^3}
       -\sum_{n=1}^N\frac1{n^2}\sum_{k=1}^{n-1}\frac{1}{n-k}
       +\sum_{n=1}^N\frac1{n^2}
        \bigg(\sum_{k=n+1}^N\frac1{k-n}-\sum_{k=1}^N\frac1k\bigg).
\end{align*}
Since $T=0$, this implies that
\[
   \sum_{n=1}^N\frac1{n^3}-\sum_{n=1}^N\frac1{n^2}\sum_{k=1}^{n-1}\frac1k
   =\sum_{n=1}^N\frac1{n^2}\bigg(\sum_{k=1}^N
   \frac1k-\sum_{k=1}^{N-n}\frac1k\bigg)
   = \sum_{n=1}^N\frac1{n^2}\sum_{k=1}^n \frac{1}{N-k+1},
\]
which is~\eqref{z21finite}.  But the right hand side satisfies
\begin{align*}
   \frac{H_N}{N}
   = \sum_{n=1}^N \frac1{n^2}\cdot\frac{n}{N}
  & \le \sum_{n=1}^N \frac1{n^2} \sum_{k=1}^n \frac1{N-k+1}\\
  & \le \sum_{n=1}^N \frac{1}{n^2}\cdot\frac{n}{N-n+1}
   =\frac1{N+1}\sum_{n=1}^N\bigg(\frac1{n}+\frac1{N-n+1}\bigg)
   =\frac{2H_N}{N+1}.
\end{align*}
Letting $N$ grow without bound now gives~\eqref{z21}, since $\dst
\lim_{N\to\infty}\frac{H_N}{N}=0$.

\eop
\section{Geometric Series}

\subsection{Convolution of Geometric Series}\label{sect:Williams} The
following argument is suggested in~\cite{Williams}.  A closely
related derivation, in which our explicit consideration of the error
term is suppressed by taking $N$ infinite, appears
in~\cite{Bracken}. Let $2\le m\in\Z$, and consider
\begin{align*}
   \sum_{j=1}^{m-2}\zeta(j+1)\zeta(m-j)
   &= \lim_{N\to\infty} \sum_{n=1}^N\sum_{k=1}^N \sum_{j=1}^{m-2}
      \frac{1}{n^{j+1}}\frac{1}{k^{m-j}}\\
   &= \lim_{N\to\infty} \bigg\{ \sum_{\substack{n,k=1\\k\ne n}}^N
      \bigg(\frac{1}{n^{m-1}(k-n)k}-\frac{1}{n(k-n)k^{m-1}}\bigg)
      +\sum_{n=1}^N\frac{m-2}{n^{m+1}}\bigg\}\\
   &=(m-2)\zeta(m+1)+2\lim_{N\to\infty}\sum_{\substack{n,k=1\\k\ne
   n}}^N\frac{1}{n^{m-1}k(k-n)}.
\end{align*}
Thus, we find that
\begin{align*}
   &   (m-2)\zeta(m+1)-\sum_{j=1}^{m-2}\zeta(j+1)\zeta(m-j)\\
   &= 2\lim_{N\to\infty}\sum_{n=1}^N\frac{1}{n^m}\sum_{\substack{k=1\\
      k\ne n}}^N \bigg(\frac1k-\frac1{k-n}\bigg)\\
   &=
   2\lim_{N\to\infty}\sum_{n=1}^N\frac1{n^m}\bigg\{\sum_{k=1}^{n-1}
   \frac1k-\frac1n+\sum_{k=1}^n\frac1{N-k+1}\bigg\}\\
   &= 2\zeta(m,1)-2\zeta(m+1)+2\lim_{N\to\infty}\sum_{n=1}^N
   \frac1{n^m}\sum_{k=1}^n\frac1{N-k+1},
\end{align*}
and hence
\[
   2\zeta(m,1)= m\zeta(m+1)-\sum_{j=1}^{m-2}\zeta(j+1)\zeta(m-j)
      -2\lim_{N\to\infty}\sum_{n=1}^N\frac1{n^m}\sum_{k=1}^n\frac1{N-k+1}.
\]
But, in light of
\[
   \sum_{n=1}^N\frac1{n^m}\sum_{k=1}^n\frac1{N-k+1}
   \le \sum_{n=1}^N \frac1{n^m}\cdot\frac{n}{N-n+1}
   \le\frac1{N+1}\sum_{n=1}^N\bigg(\frac1{N-n+1}+\frac1n\bigg)
   =\frac{2H_N}{N+1},
\]
the identity~\eqref{EulerReduction} now follows. \qed

\subsection{A Sum Formula}\label{sect:SumFormula}
Equation~\eqref{z21} is the case $n=3$ of the following result.
See~\cite{Briggs0}.
\begin{Thm} \label{briggs} If $3\le n\in\Z$ then
\begin{equation}
   \zeta(n) = \sum_{j=1}^{n-2}\zeta(n-j,j).
\label{sumdepth2}
\end{equation}
\end{Thm}

We discuss a generalization~\eqref{sum} of the sum
formula~\eqref{sumdepth2} to arbitrary depth in
\S\ref{sect:SumGF}.

\noindent{\bf Proof.}  Summing the geometric series on the right
hand side gives
\begin{align*}
   \sum_{j=1}^{n-2}\sum_{h=1}^\infty\sum_{m=1}^\infty
   \frac1{h^j(h+m)^{n-j}}
   &= \sum_{h,m=1}^\infty
   \bigg[\frac{1}{h^{n-2}m(h+m)}-\frac{1}{m(h+m)^{n-1}}\bigg]\\
   &=\sum_{h=1}^\infty \frac1{h^{n-1}}\sum_{m=1}^\infty \bigg(
   \frac1m-\frac1{h+m}\bigg)-\zeta(n-1,1)\\
   &= \sum_{h=1}^\infty \frac1{h^{n-1}}\sum_{k=1}^{h}\frac1k
   -\zeta(n-1,1)\\
   &= \sum_{h=1}^\infty \frac1{h^n} +\sum_{h=1}^\infty
   \frac1{h^{n-1}}\sum_{k=1}^{n-1}\frac1k-\zeta(n-1,1)\\
   &= \zeta(n).
\end{align*} \qed

\subsection{A $q$-Analogue}\label{sect:zud} The
following argument is based on an idea of Zudilin~\cite{Zud}. We
begin with the finite geometric series identity
\[
   \frac{uv}{(1-u)(1-uv)^s}+\frac{uv^2}{(1-v)(1-uv)^s}
   = \frac{uv}{(1-u)(1-v)^s} -
   \sum_{j=1}^{s-1}\frac{uv^2}{(1-v)^{j+1}(1-uv)^{s-j}},
\]
valid for all positive integers $s$ and real $u$, $v$ with $u\ne
1$, $uv\ne 1$.  We now assume $s>1$, $q$ is real and $0<q<1$. Put
$u=q^m$, $v=q^n$ and sum over all positive integers $m$ and $n$.
Thus,
\begin{align*}
   &\sum_{m,n>0} \frac{q^{m+n}}{(1-q^m)(1-q^{m+n})^s}
   +\sum_{m,n>0}\frac{q^{m+2n}}{(1-q^n)(1-q^{m+n})^s}\\
   &= \sum_{m,n>0}\frac{q^{m+n}}{(1-q^m)(1-q^n)^s}
   -\sum_{m,n>0}\frac{q^{m+2n}}{(1-q^n)^s(1-q^{m+n})}\\
   &\qquad -\sum_{j=1}^{s-2}\sum_{m,n>0}
   \frac{q^{m+2n}}{(1-q^n)^{j+1}(1-q^{m+n})^{s-j}}\\
   &= \sum_{m,n>0}\frac{q^n}{(1-q^n)^s}\bigg[\frac{q^m}{1-q^m}-
   \frac{q^{m+n}}{1-q^{m+n}}\bigg]-\sum_{j=1}^{s-2}\sum_{m,n>0}
   \frac{q^{m+2n}}{(1-q^n)^{j+1}(1-q^{m+n})^{s-j}}\\
   &= \sum_{n>0}\frac{q^n}{(1-q^n)^s}\sum_{m=1}^n\frac{q^m}{1-q^m}
   -\sum_{j=1}^{s-2}\sum_{m,n>0}
   \frac{q^{m+2n}}{(1-q^n)^{j+1}(1-q^{m+n})^{s-j}}\\
   &=\sum_{n>0}\frac{q^{2n}}{(1-q^n)^{s+1}}+\sum_{n>m>0}
   \frac{q^{n+m}}{(1-q^n)^s(1-q^m)}-\sum_{j=1}^{s-2}\sum_{m,n>0}
   \frac{q^{m+2n}}{(1-q^n)^{j+1}(1-q^{m+n})^{s-j}}\\
\end{align*}
Cancelling the second double sum on the left with the
corresponding double sum on the right and replacing $m+n$ by $k$
in the remaining sums now yields
\[
   \sum_{k>m>0}\frac{q^{k}}{(1-q^{k})^s(1-q^m)}
   = \sum_{n>0}\frac{q^{2n}}{(1-q^n)^{s+1}}-
   \sum_{j=1}^{s-2}\sum_{k>m>0}
   \frac{q^{k+m}}{(1-q^m)^{j+1}(1-q^{k})^{s-j}},
\]
or equivalently, that
\begin{equation}
   \sum_{k>0}\frac{q^{2k}}{(1-q^k)^{s+1}}
   = \sum_{k>m>0}\frac{q^k}{(1-q^k)^s(1-q^m)}
   +\sum_{j=1}^{s-2}\sum_{k>m>0}\frac{q^{k+m}}{(1-q^k)^{s-j}(1-q^m)^{j+1}}.
\label{general}
\end{equation}
Multiplying~\eqref{general} through by $(1-q)^{s+1}$ and letting
$q\to 1$ gives
\[
   \zeta(s+1)=\zeta(s,1)+\sum_{j=1}^{s-2} \zeta(s-j,j+1),
\]
which is just a restatement of~\eqref{sumdepth2}. Taking $s=2$
gives~\eqref{z21} again. \qed

As in~\cite{DBq}, define the $q$-analog of a non-negative integer
$n$ by
\[
   [n]_q := \sum_{k=0}^{n-1} q^k = \frac{1-q^n}{1-q},
\]
and the multiple $q$-zeta function
\begin{equation}
   \zeta[s_1,\dots,s_m] := \sum_{k_1>\cdots >k_m>0}
   \; \prod_{j=1}^m \frac{q^{(s_j-1)k_j}}{[k_j]_q^{s_j}},
   \label{qMZVdef}
\end{equation}
where $s_1,s_2,\dots,s_m$ are real numbers with $s_1>1$ and
$s_j\ge 1$ for $2\le j\le m$.  Then multiplying~\eqref{general} by
$(1-q)^{s+1}$ and then setting $s=2$ gives $\zeta[2,1]=\zeta[3]$,
which is a $q$-analog of~\eqref{z21}.  That is, the latter may be
obtained from the former by letting $q\to 1-$.  On the other hand,
$s=3$ in~\eqref{general} gives
\[
   \zeta[4] +(1-q)\zeta[3]=\zeta[3,1]+(1-q)\zeta[2,1]+\zeta[2,2],
\]
which, in light of $\zeta[2,1]=\zeta[3]$ implies $\zeta[3,1] =
\zeta[4]-\zeta[2,2]$.  By Theorem 1 of~\cite{DBq}, we know that
$\zeta[2,2]$ reduces to depth 1 multiple $q$-zeta values.  Indeed,
by the $q$-stuffle multiplication rule~\cite{DBq},
$
   \zeta[2]\zeta[2] = 2\zeta[2,2] + \zeta[4] +(1-q)\zeta[3].
$
Thus,
\[
   \zeta[3,1] = \zeta[4]-\zeta[2,2]=
   \tfrac32\zeta[4]-\tfrac12\left(\zeta[2]\right)^2+\tfrac12(1-q)\zeta[3],
\]
which is a $q$-analog of the
evaluation~\cite{BBB,BBBLa,BBBLc,BowBrad1,BowBradSurvey,BowBrad3,BowBradRyoo}
\[
   \zeta(3,1) = \frac{\pi^4}{360}.
\]
Additional material concerning $q$-analogs of multiple harmonic sums
and multiple zeta values can be found
in~\cite{DBq,DBqKarl,DBqSum,DBqDecomp}.

\section{Integral Representations}\label{sect:integrals}
\subsection{Single Integrals I}\label{sect:single1}

We use the fact that
\begin{equation}\label{naive}
   \int_0^1 u^{k-1}(-\log u)\,du = \frac1{k^2},\qquad k>0.
\end{equation}
Thus
\begin{align}\label{logs}
   \sum_{k>n>1}\frac{1}{k^2n}
   &= \sum_{n=1}^\infty \frac1n\sum_{k>n}\int_0^1 u^{k-1}(-\log u)\,du
      \nonumber\\
   &= \sum_{n=1}^\infty\frac1n \int_0^1 (-\log u)\sum_{k>n} u^{k-1}\,du
       \nonumber\\
   &= \sum_{n=1}^\infty\frac1n \int_0^1 (-\log u)
                                \frac{u^n}{1-u}\,du \nonumber\\
   &= -\int_0^1 \frac{\log u}{1-u}\sum_{n=1}^\infty \frac{u^n}{n}\,du
       \nonumber\\
   &= \int_0^1 (-\log u)(1-u)^{-1}\log(1-u)^{-1}\,du.
\end{align}
The interchanges of summation and integration are in each case
justified by Lebesgue's monotone convergence theorem. After making
the change of variable $t=1-u$, we obtain
\begin{equation}\label{morelogs}
  \sum_{k>n>1}\frac{1}{k^2n}
   = \int_0^1 \log(1-t)^{-1} (-\log t)\,\frac{dt}{t}
   = \int_0^1 (-\log t) \sum_{n=1}^\infty \frac{t^{n-1}}{n}\,dt.
\end{equation}
Again, since all terms of the series are positive, Lebesgue's
monotone convergence theorem permits us to interchange the order
of summation and integration.  Thus, invoking~\eqref{naive} again,
we obtain
\[
   \sum_{k>n>1}\frac{1}{k^2n}
   = \sum_{n=1}^\infty \frac1n\int_0^1 (-\log t)\, t^{n-1}\,dt
   = \sum_{n=1}^\infty \frac1{n^3},
\]
which is~\eqref{z21}. \qed

\subsection{Single Integrals II}\label{sect:single2}
The Laplace transform
\begin{equation}
   \int_0^1 x^{r-1} (-\log x)^\sigma\,dx
   = \int_0^{\infty} e^{-ru}\, u^\sigma\,du
   = \frac{\Gamma(\sigma+1)}{r^{\sigma+1}},\qquad r>0,\quad \sigma>-1,
\label{Laplace}
\end{equation}
generalizes~\eqref{naive} and yields the representation
\[
   \zeta(m+1) =
   \frac{1}{m!}\sum_{r=1}^\infty\frac{\Gamma(m+1)}{r^{m+1}}
   = \frac{1}{m!}\sum_{r=1}^\infty \int_0^1 x^{r-1}(-\log
   x)^m \,dx
   = \frac{(-1)^m}{m!}\int_0^1 \frac{\log^m x}{1-x}\,dx.
\]
The interchange of summation and integration is valid if $m>0$.
The change of variable $x\mapsto 1-x$ now yields
\begin{equation}
   \zeta(m+1) = \frac{(-1)^m}{m!}\int_0^1
   \log^m(1-x)\,\frac{dx}{x},
   \qquad 1\le m\in\Z.
\label{ZetaStirling}
\end{equation}

In~\cite{Farnum}, equation~\eqref{Laplace} in conjunction with
clever use of change of variable and integration by parts, is used
to prove the identity
\begin{equation}
   k!\zeta(k+2) = \sum_{n_1=1}^\infty\, \sum_{n_2=1}^\infty \cdots
   \sum_{n_k=1}^\infty \frac{1}{n_1n_2\cdots n_k}\;
   \sum_{p=1+n_1+n_2+\cdots+n_k}\frac1{p^2},
   \qquad 0\le k\in\Z.
\label{Tissier}
\end{equation}
The case $k=1$ of~\eqref{Tissier} is precisely~\eqref{z21}.  We
give here a slightly simpler proof of~\eqref{Tissier}, dispensing
with the integration by parts.

From~\eqref{Laplace},
\begin{align*}
   k!\zeta(k+2) &= \sum_{r=1}^\infty
   \frac1{r}\cdot\frac{\Gamma(k+1)}{r^{k+1}}
   = \sum_{r=1}^\infty \frac1r\int_0^1 x^{r-1}(-\log x)^k\,dx\\
   &= \int_0^1 (-\log x)^k \log(1-x)^{-1}\,\frac{dx}{x}\\
   &= \int_0^1 \log^k(1-x)^{-1} (-\log x)\frac{dx}{1-x}\\
   &=\sum_{n_1=1}^\infty\, \sum_{n_2=1}^\infty \cdots\sum_{n_k=1}^\infty \frac{1}{n_1n_2\cdots n_k}\;
     \int_0^1 \frac{x^{n_1+n_2+\cdots+n_k}}{1-x}(-\log x)\,dx\\
   &= \sum_{n_1=1}^\infty\, \sum_{n_2=1}^\infty \cdots\sum_{n_k=1}^\infty \frac{1}{n_1n_2\cdots n_k}\;
     \sum_{p>n_1+n_2+\cdots+n_k} \int_0^1 x^{p-1}(-\log x)\,dx\\
   &= \sum_{n_1=1}^\infty\, \sum_{n_2=1}^\infty \cdots\sum_{n_k=1}^\infty \frac{1}{n_1n_2\cdots n_k}\;
     \sum_{p>n_1+n_2+\cdots+n_k}\frac1{p^2}.
\end{align*}
\qed

\subsection{Double Integrals I}\label{sect:double1}
Write
\begin{align*}
   \zeta(2,1) =\sum_{k,m>0}\frac{1}{k(m+k)^2}
   &= \int_0^1\int_0^1 \sum_{k>0}\frac{(xy)^k}{k}\sum_{m>0}
   (xy)^{m-1}\,dx\,dy\\
   &= -\int_0^1\int_0^1 \frac{\log(1-xy)}{1-xy}\,dx\,dy.
\end{align*}
Now make the change of variable $u=xy$, $v=x/y$ with Jacobian
$1/(2v)$, obtaining
\begin{align*}
   \zeta(2,1) = -\frac12\int_0^1
   \frac{\log(1-u)}{1-u}\int_u^{1/u}\frac{dv}{v}\,du
   = \int_0^1 \frac{(\log u)\log(1-u)}{1-u}\,du,
\end{align*}
which is~\eqref{logs}.  Now continue as in \S\ref{sect:single1}.
\qed

\subsection{Double Integrals II}\label{sect:double2}
The following is reconstructed from a phone conversation with
Krishna Alladi.  See also~\cite{Beukers}.  Let $\e>0$. By
expanding the integrand as a geometric series, one sees that
\[
   \sum_{n=1}^\infty \frac1{(n+\e)^2} = \int_0^1\int_0^1
   \frac{(xy)^{\e}}{1-xy}\,dx\,dy.
\]
Differentiating with respect to $\e$ and then letting $\e=0$ gives
\[
   \zeta(3) = -\frac12\int_0^1\int_0^1 \frac{\log(xy)}{1-xy}\,dx\,dy
   = -\frac12\int_0^1\int_0^1 \frac{\log x+\log y}{1-xy}\,dx\,dy
   = -\int_0^1\int_0^1\frac{\log x}{1-xy}\,dx
\]
by symmetry.  Now integrate with respect to $y$ to get
\begin{equation}\label{parts}
   \zeta(3) = \int_0^1 (\log x)\log(1-x)\frac{dx}{x}.
\end{equation}
Comparing~\eqref{parts} with~\eqref{morelogs} completes the proof
of~\eqref{z21}. \qed

\subsection{Integration by Parts}\label{sect:parts}
Start with~\eqref{parts} and integrate by parts, obtaining
\[
   2\zeta(3) = \int_0^1\frac{\log^2 x}{1-x}\,dx
   = \int_0^1 \log^2(1-x)\frac{dx}{x}
   = \sum_{n,k>0} \int_0^1 \frac{x^{n+k-1}}{nk}\,dx
   = \sum_{n,k>0}\frac{1}{nk(n+k)}.
\]
Now see \S\ref{sect:telescope}.  \qed

\subsection{Triple Integrals I}\label{sect:iterint1}
This time, instead of~\eqref{naive} we use the elementary identity
\[
   \frac{1}{k^2n} = \int_0^1 y_1^{-1}\int_0^{y_1}y_2^{k-n-1}
   \int_0^{y_2} y_3^{n-1}\,dy_3\,dy_2\,dy_1,\qquad k>n>0.
\]
This yields
\begin{equation}\label{irep}
    \sum_{k>n>0}\frac{1}{k^2n}
   = \int_0^1 y_1^{-1}\int_0^{y_1}(1-y_2)^{-1}\int_0^{y_2}(1-y_3)^{-1}
      \,dy_3\,dy_2\,dy_1.
\end{equation}
Now make the change of variable $y_i=1-x_i$ for $i=1,2,3$ to
obtain
\begin{align*}
   \sum_{k>n>0}\frac{1}{k^2n}
   &= \int_0^1 (1-x_1)^{-1}\int_{x_1}^1 x_2^{-1}\int_{x_2}^1 x_3^{-1}
      \,dx_3\,dx_2\,dx_1\\
   &= \int_0^1 x_3^{-1}\int_0^{x_3}x_2^{-1}\int_0^{x_2}(1-x_1)^{-1}
      \,dx_1\,dx_2\,dx_3.
\end{align*}
After expanding $(1-x_1)^{-1}$ into a geometric series and
interchanging the order of summation and integration, one arrives
at
\[
   \sum_{k>n>0}\frac{1}{k^2n}
   = \sum_{n=1}^\infty \int_0^1 x_3^{-1}\int_0^{x_3}x_2^{-1}\int_0^{x_2}
       x_1^{n-1}\,dx_1\,dx_2\,dx_3
   = \sum_{n=1}^\infty \frac1{n^3},
\]
as required. \qed

More generally~\cite{BBB,BBBLa,BBBLc,BowBradSurvey,BowBrad3,CK},
\begin{equation}
   \zeta(s_1,\dots,s_k)
   = \sum_{n_1>\cdots>n_k>0}\;\prod_{j=1}^k n_j^{-s_j}
   =\int \prod_{j=1}^k
   \bigg(\prod_{r=1}^{s_j-1} \frac{dt_r^{(j)}}{t_r^{(j)}}\bigg)
   \frac{dt_{s_j}^{(j)}}{1-t_{s_j}^{(j)}},
\label{iterint1}
\end{equation}
where the integral is over the simplex
\[
   1>t_1^{(1)}>\cdots>t_{s_1}^{(1)}>\cdots>t_1^{(k)}>\cdots>t_{s_k}^{(k)}>0,
\]
and is abbreviated by
\begin{equation}
   \int_0^1 \prod_{j=1}^k a^{s_j-1}b,
   \qquad a=\frac{dt}{t},\qquad b=\frac{dt}{1-t}.
\label{shortiterint1}
\end{equation}
The change of variable $t\mapsto 1-t$ at each level of integration
switches the differential forms $a$ and $b$, thus yielding the
duality formula~\cite{BBB}~\cite[p.\ 483]{CK} (conjectured
in~\cite{Hoff92})
\begin{equation}\label{duality}
   \zeta(s_1+2,\us^{r_1},\dots,s_n+2,\us^{r_n})
      = \zeta(r_n+2,\us^{s_n},\dots,r_1+2,\us^{s_1}),
\end{equation}
which is valid for all nonnegative integers
$s_1,r_1,\ldots,s_n,r_n$.  The case $s_1=0$, $r_1=1$
of~\eqref{duality} is~\eqref{z21}.  More generally,~\eqref{z21^n}
can be restated as
\[
   \int_0^1 (ab^2)^n = \int_0^1 (a^2b)^n
\]
and thus~\eqref{z21^n} is recovered by taking each $s_j=0$ and
each $r_j=1$ in~\eqref{duality}.  For further generalizations and
extensions of duality, see~\cite{BBBLa,DBq,DBqKarl}.

For alternations, we require in addition the differential form
$c:=-dt/(1+t)$ with which we may form the generating function
\[
   \sum_{n=1}^\infty z^{3n} \zeta(\{\overline2,1\}^n)
   = \sum_{n=0}^\infty \bigg\{ z^{6n+3}\int_0^1
   (ac^2ab^2)^n ac^2 +  z^{6n+6} \int_0^1
   (ac^2ab^2)^{6n+6}\bigg\}.
\]
A lengthy calculation verifies that the only changes of variable
that preserve the unit interval and send the non-commutative
polynomial ring $\Q\langle a,b\rangle$ into $\Q\langle
a,b,c\rangle$ are
\begin{alignat}{2}
   S(a,b) &= S(a,b),\qquad\qquad & t &\mapsto t,\label{identity}\\
   S(a,b) &= R(b,a),\qquad\qquad & t &\mapsto 1-t,\label{tau}\\
   S(a,b) &= S(2a,b+c),\qquad\qquad & t &\mapsto t^2,\label{sumsigns}\\
   S(a,b) &= S(a+c,b-c),\qquad\qquad & t &\mapsto
               \frac{2t}{1+t},\label{Landen}\\
   S(a,b) &= S(a+2c,2b-2c),\qquad\qquad & t &\mapsto
             \frac{4t}{(1+t)^2},\label{quadLanden}
\end{alignat}
and compositions thereof, such as $t\mapsto
1-2t/(1+t)=(1-t)/(1+t)$, etc.
In~\eqref{identity}--\eqref{quadLanden}, $S(a,b)$ denotes a
non-commutative word on the alphabet $\{a,b\}$ and $R(b,a)$
denotes the word formed by switching $a$ and $b$ and then
reversing the order of the letters.

Now view $a$, $b$ and $c$ as indeterminates.  In light of the
polynomial \emph{identity}
\[
   ab^2-8ac^2 = 2[ab^2-2a(b+c)^2] + 8[ab^2-(a+c)(b-c)^2]
   + [(a+2c)(2b-2c)^2-ab^2]
\]
in the non-commutative ring $\Z\langle a,b,c\rangle$  and the
transformations~\eqref{sumsigns}, \eqref{Landen} and
\eqref{quadLanden} above, each bracketed term vanishes when we make
the identifications $a=dt/t$, $b=dt/(1-t)$, $c=-dt/(1+t)$ and
perform the requisite iterated integrations.   Thus,
\[
   \zeta(2,1)-8\zeta(\overline2,1)=\int_0^1 ab^2-8\int_0^1 ac^2 =0,
\]
which in light of~\eqref{z21} proves~\eqref{2bar1}.  \qed

\subsection{Triple Integrals II}\label{sect:iterint2}
First, note that by expanding the integrands in geometric series
and integrating term by term,
\[
   \zeta(2,1) = 8\int_0^1 \frac{dx}{x}\int_0^x
   \frac{y\,dy}{1-y^2}\int_0^y\frac{z\,dz}{1-z^2}.
\]
Now make the change of variable
\[
   \frac{x\,dx}{1-x^2}=\frac{du}{1+u},
   \qquad
   \frac{y\,dy}{1-y^2}=\frac{dv}{1+v},
   \qquad
   \frac{z\,dz}{1-z^2}=\frac{dw}{1+w}
\]
to obtain the equivalent integral
\[
   \zeta(2,1) = 8\int_0^\infty
   \bigg(\frac{du}{2u}+\frac{du}{2(2+u)}
   -\frac{du}{1+u}\bigg)\int_0^u\frac{dv}{1+v}
   \int_0^v\frac{dw}{1+w}.
\]
The two inner integrals can be directly performed, leading to
\[
   \zeta(2,1) = 4\int_0^\infty\frac{\log^2(u+1)}{u(u+1)(u+2)}\,{du}.
\]
Finally, make the substitution $u+1=1/\sqrt{1-x}$ to obtain
\[
   \zeta(2,1) = \frac12\int_0^1 \frac{\log^2(1-x)}{x}\,dx
   = \zeta(3),
\]
by~\eqref{ZetaStirling}.

\subsection{Complex Line Integrals I}\label{sect:Perron}
Here we apply the Mellin inversion formula~\cite[p.\
243]{Apostol}, \cite[pp.\ 130--132 ]{Tenenbaum}
\[
   \frac{1}{2\pi i}\int_{c-i\infty}^{c+i\infty} y^z\,\frac{dz}{z} =
\   \begin{cases} 1,\quad y>1 \\ 0, \quad y<1\\ \tfrac12,\quad
   y=1\end{cases}
\]
which is valid for fixed $c>0$. It follows that if $c>0$ and
$s-1>c>1-t$ then the Perron-type formula
\begin{align}
   \zeta(s,t) +\frac12\zeta(s+t)
   &= \sum_{n=1}^\infty n^{-s} \sum_{k=1}^\infty k^{-t}
   \frac1{2\pi i} \int_{c-i\infty}^{c+i\infty}
   \bigg(\frac{n}{k}\bigg)^z\frac{dz}{z}\nonumber\\
   &=\frac1{2\pi i}\int_{c-i\infty}^{c+i\infty}
   \zeta(s-z)\zeta(t+z)\,\frac{dz}{z}
\label{Perron}
\end{align}
is valid.  (Interchanging the order of summation and integration
is permissible by absolute convergence.)  Although we have not yet
found a way to exploit~\eqref{Perron} in proving identities such
as~\eqref{z21}, we note that by integrating around the rectangular
contour with corners $(\pm c\pm iM)$ and then letting
$M\to+\infty$, one can readily establish the
stuffle~\cite{BBBLa,BowBradSurvey,DBPrtn,DBq} formula in the form
\[
   \zeta(s,t)+\frac12\zeta(s+t)+\zeta(t,s)+\frac12\zeta(t+s)
   = \zeta(s)\zeta(t),  \qquad s,t>1+c.
\]
The right hand side arises as the residue contribution of the
integrand at $z=0$.  One can also use~\eqref{Perron} to establish
\[
   \sum_{s=2}^\infty\left[\zeta(s,1)+\tfrac12\zeta(s+1)\right]x^{s-1}
   = \sum_{n>m>0}\frac{x}{mn(n-x)} +\frac12\sum_{n=1}^\infty
   \frac{x}{n(n-x)},
\]
but this is easy to prove directly.

\subsection{Complex Line Integrals II}\label{sect:Dirichlet}
We let $\lambda(s) :=\sum_{n>0} \lambda_n\, n^{-s}$ represent a
formal Dirichlet series, with real coefficients $\lambda_n$ and we
set $s:=\sigma+i\,\tau$ with $\sigma=\Re(s)>0$, and consider the
following integral:
\begin{eqnarray}\label{int1}
\iota_\lambda(\sigma):=\int_{0}^\infty
\left|\frac{\lambda(s)}s\right|^2\,d\tau = \frac 12\,
\int_{-\infty}^\infty \left|\frac{\lambda(s)}s\right|^2\,d\tau,
\end{eqnarray}
as a function of $\lambda$.  We begin with a useful variant of the
Mellin inversion formula, namely
\begin{eqnarray}\label{intc}\int _{-\infty}^{\infty }\!{\frac
{\cos \left( at \right) }{{t}^{2}+{ u}^{2}}}{dt}={\frac {\pi
}{u}}\,e^{-au},
\end{eqnarray}
for $u, a>0$, as follows by contour integration, from a computer
algebra system, or otherwise. This leads to

\begin{Thm}\label{dir-thm} {\rm (Theorem 1 of  \cite{JB}).}  {\it  For
$\lambda(s)=\sum_{n=1}^\infty \,\lambda_n\,n^{-s}$ and
$s=\sigma+i\,\tau$ with fixed $\sigma=\Re(s)>0$ such that the
Dirichlet series is absolutely convergent it is true that
\begin{eqnarray}\label{ans1}
\iota_\lambda(\sigma) = \int_{0}^\infty
\left|\frac{\lambda(s)}s\right|^2\,d\tau
=\frac{\pi}{2\sigma}\,\sum_{n=1}^\infty\frac{\Lambda_n^2-
\Lambda_{n-1}^2}{n^{2\sigma}},\end{eqnarray}
 where
$\Lambda_n:=\sum_{k=1}^{n}\lambda_k$ and $\Lambda_0:=0$.}

\emph{More generally, for given absolutely convergent Dirichlet
series $\alpha(s):=\sum_{n=1}^\infty\,{\alpha_n}\,{n^{-s}}$ and
$\beta(s):=\sum_{n=1}^\infty\,{\beta_n}\,{n^{-s}}$
\begin{eqnarray}\label{ans-gen}  \frac12\int_{-\infty}^\infty
\frac{\alpha(s)\,\overline{\beta}(s)}{\sigma^2+\tau^2}\,d\tau =
\frac{\pi}{2\sigma}\,\sum_{n=1}^\infty
\frac{A_n\,\overline{B_n}-A_{n-1}\,\overline{B_{n-1}}}{n^{2\sigma}},\label{ans3}
\end{eqnarray}
in which $A_n=\sum_{k=1}^{n}\alpha_k$ and
$B_n=\sum_{k=1}^{n}\beta_k.$}
 \end{Thm}

Note that the righthand side of \eqref{ans1} is always a generalized
Euler sum. \qed

 For the Riemann zeta function, and for
$\sigma>1$, Theorem \ref{dir-thm} applies and yields
    $$
\frac{\sigma}{\pi}\,\iota_\zeta(\sigma)=
    \zeta(2\sigma-1)-\frac12\,\zeta(2\sigma),$$ as $\lambda_n=1$ and
    $\Lambda_n=n-1/2$.  By contrast it is known that on the critical
    line
$$\frac{1/2}{\pi}\,\iota_\zeta\left(\frac 12\right)=\log(\sqrt{2\,\pi})-\frac 12\,\gamma.$$
 There are similar formulae for
    $s \mapsto \zeta(s-k)$ with $k$ integral. For instance,
    applying the result in (\ref{ans1}) with $\zeta_1:=t \mapsto \zeta(t+1)$ yields
    $$\frac{1}{\pi}\,\int_{0}^\infty
\frac{|\zeta(3/2+i\tau)|^2}{1/4+\tau^2}\,d\tau =
\frac{1}{\pi}\,\iota_{\zeta_1}\left(\frac 1 2\right)=
 2\, \zeta(2,1)+\zeta(3)=3\,\zeta(3),$$
   on using ~\eqref{z21}. For the \emph{alternating zeta function},
    $\alpha:=s\mapsto(1-2^{1-s})\zeta(s)$,
     the same approach  via (\ref{ans-gen}) produces
         $$\frac{1}{\pi}\,\int_{0}^\infty
\frac{\alpha(3/2+i\tau)\,\overline{\alpha(3/2+i\tau)}}{1/4+\tau^2}\,d\tau
=
 2\,\zeta(\overline{2},\overline{1})+\zeta(3)=3\, \zeta(2)\,\log(2)-\frac 94 \zeta(3),$$
 and
      $$\frac{1}{2\pi}\,\int_{-\infty}^\infty
\frac{\alpha(3/2+i\tau)\,\overline{\zeta(3/2+i\tau)}}{1/4+\tau^2}\,d\tau
=
 \zeta(\overline{2},1)+\zeta(2,\overline{1})+\alpha(3)=\frac 98\, \zeta(2)\,\log(2)-\frac {3}4 \zeta(3),$$
 since as we have seen repeatedly $\zeta(\overline{2},1)
    =\zeta(3)/8$; while $\zeta(2,\overline{1})=\zeta(3)-3/2\,\zeta(2)\log(2)$ and  $\zeta(\overline{2},\overline{1})=
3/2\,\zeta(2)\log(2)-13/8\,\zeta(3),$ (e.g., \cite{BZB}).

As in the previous subsection we have not been able to directly
obtain \eqref{2bar1} or even ~\eqref{z21}, but we have connected
them to quite difficult line integrals.

\subsection{Contour Integrals and Residues}\label{sect:residue}
Following~\cite{State}, let $\mathscr{C}_n$ $(n\in\Z^{+})$ be the
square contour with vertices $(\pm 1\pm i)(n+1/2)$.  Using the
asymptotic expansion
\[
   \psi(z) \sim \log z -
   \frac1{2z}-\sum_{r=1}^\infty\frac{B_{2r}}{2rz^{2r}},
   \qquad |\arg z|<\pi
\]
in terms of the Bernoulli numbers
\[
   \frac{t}{1-e^{-t}}=1+\frac{t}{2}+\sum_{r=1}^\infty
   \frac{B_{2r}}{(2r)!}t^{2r},\qquad |t|<2\pi
\]
and the identity
\[
   \psi(z)=\psi(-z) - \frac1{z}-\pi\cot \pi z,
\]
we can show that for each integer $k\ge 2$,
\[
   \lim_{n\to\infty} \int_{\mathscr{C}_n} z^{-k}\,\psi^2(-z)\,dz =
   0.
\]
Then by the residue theorem, we obtain
\begin{Thm}[Theorem 3 of~\cite{State}]\label{thm:State} For every integer $k\ge
2$,
\[
   2\sum_{n=1}^\infty n^{-k}\,\psi(n) = k\zeta(k+1)-2\gamma
   \zeta(k)-\sum_{j=1}^{k-1}\zeta(j)\zeta(k-j+1),
\]
where $\gamma=0.577215664\dots$ is Euler's constant.
\end{Thm}
In light of the identity
\[
   \psi(n)+\gamma= H_{n-1}=\sum_{k=1}^{n-1}\frac1k,\qquad
   n\in\Z^{+},
\]
Theorem~\ref{thm:State} is equivalent to~\eqref{EulerReduction}.
The case $k=2$ thus gives~\eqref{z21}. \qed

Flajolet and Salvy~\cite{Flaj} developed the residue approach more
systematically, and applied it to a number of other Euler sum
identities in addition to~\eqref{EulerReduction}.

\section{Witten Zeta-functions}\label{sect:witten} We recall that  for
$r,s>1/2$:
$$\mathcal{W}(r,s,t):= \sum_{n=1}^\infty\sum_{m=1}^\infty \frac{1}{n^r\,m^s\,(n+m)^t}$$
is a \emph{Witten $\zeta$-function}, \cite{zagier,moll,CranBuh}. We
refer to \cite{zagier} for a description of the uses of more general
Witten $\zeta$-functions. Ours are also called \emph{Tornheim double
sums}, \cite{moll}. There is a simple algebraic relation
\begin{eqnarray}\label{w-alg}\mathcal{W}(r,s,t)=\mathcal{W}(r-1,s,t+1)+\mathcal{W}(r,s-1,t+1).\end{eqnarray}
This is based on writing
$$\frac{m+n}{(m+n)^{t+1}}=\frac{m}{(m+n)^{t+1}}+\frac{n}{(m+n)^{t+1}}.$$ Also
\begin{eqnarray}\label{w-alg1}\mathcal{W}(r,s,t) =
\mathcal{W}(s,r,t),\end{eqnarray} and
\begin{eqnarray}\label{w-alg2}\mathcal{W}(r,s,0)=\zeta(r)\,\zeta(s)\quad
\mbox{while} \quad \mathcal{W}(r,0,t) =\zeta(t,r).\end{eqnarray}

Hence, $\mathcal{W}(s,s,t)=2\,\mathcal{W}(s,s-1,t+1)$ and so
$$\mathcal{W}(1,1,1)=2\,\mathcal{W}(1,0,2)=2\,\zeta(2,1)=2\,\zeta(3).$$
Note  the analogue to (\ref{w-alg}), viz.
$\zeta(s,t)+\zeta(t,s)=\zeta(s)\,\zeta(t)-\zeta(s+t)$, shows
$\mathcal{W}(s,0,s)=2\,\zeta(s,s)=\zeta^2(s)-\zeta(2s).$ Thus,
$\mathcal{W}(2,0,2)=2\,\zeta(2,2)=\pi^4/36-\pi^4/90=\pi^4/72$.

More generally, recursive use of (\ref{w-alg}) and (\ref{w-alg1}),
along with initial conditions (\ref{w-alg2}) shows that  \emph{all
integer $\mathcal{W}(s,r,t)$ values are expressible in terms of
double (and single) Euler sums.} If we start with
$\Gamma(s)/(m+n)^{t} = \int _0^1\! (-\log
\sigma)^{t-1}\,\sigma^{m+n-1}\,d \sigma $ we obtain
\begin{eqnarray}\label{gamma}\mathcal{W}(r,s,t)=
 \frac{1}{\Gamma(t)}\,\int _0^1\!
 {\rm Li}_r(\sigma)\,{\rm Li}_s(\sigma)\,\frac{\left(-\log \sigma\right)^{t-1}}{\sigma}\,d \sigma.\end{eqnarray}
For example, we recover an analytic proof of
\begin{eqnarray}\label{z21-w}2\,\zeta(2,1)=\mathcal{W}(1,1,1)=
 \int _0^1\!\frac{\ln^2(1-\sigma)}{\sigma}\,d\sigma =2\,\zeta(3),\end{eqnarray}
 Indeed $S$ in the proof of \S \ref{sect:telescope} is precisely $\mathcal{W}(1,1,1)$.

 We may now discover analytic as opposed to algebraic
 relations. Integration by parts yields
\begin{eqnarray}\label{w-parts}\mathcal{W}(r,s+1,1) +\mathcal{W}(r+1,s,1)={\rm Li}_{r+1}(1)\,{\rm
 Li}_{s+1}(1)=\zeta(r+1)\,\zeta(s+1),\end{eqnarray}
 So, in particular, $\mathcal{W}(s+1,s,1)=\zeta^2(s+1)/2$.

  Symbolically, \emph{Maple} immediately evaluates
$\mathcal{W}(2,1,1)=\pi^4/72,$ and while it fails directly with
$\mathcal{W}(1,1,2)$, we know it must be a multiple of $\pi^4$ or
equivalently $\zeta(4)$; and numerically obtain
$\mathcal{W}(1,1,2)/\zeta(4)=.49999999999999999998\ldots$.

 \subsection{The Hilbert Matrix} Letting $a_n:=1/n^r$ and $b_n:=1/n^s$,
 inequality (\ref{hilbert-p}) of Section \ref{sect:hardy} yields
\begin{eqnarray}\label{wittenp}\mathcal{W}(r,s,1) \le
 \pi\,{\rm csc}\left(\frac{\pi}p\right)\,\sqrt[p]{\zeta(pr)}\,\sqrt[q]{\zeta(qs)}.\end{eqnarray}
Indeed, the constant in (\ref{hilbert-p}) is best possible
\cite{ghh, steele}. We consider
$$\mathcal{R}_p(s):=\frac{\mathcal{W}((p-1)s,s,1)}{\pi\,\zeta(ps)},$$
and observe that with $\sigma_n^p(s):=\sum_{m=1}^\infty
(n/m)^{-(p-1)s}/(n+m) \to \pi\,{\rm csc}\left(\frac{\pi}q \right),$
we have
\begin{eqnarray*}\mathcal{L}_p:&=&\lim_{s\to
1/p}(ps-1)\,\sum_{n=1}^\infty\sum_{m=1}^\infty
\frac{n^{-s}\,m^{-(p-1)s}}{n+m}= \lim_{s\to
1/p}(ps-1)\,\sum_{n=1}^\infty\frac{1}{n^{ps}}\,\sigma_n^p(s)\\&=&\lim_{s\to
1/p}\,(ps-1) \,\sum_{n=1}^\infty
\,\frac{\left\{\sigma_n^p(s)-\pi\,{\rm csc}\left(\pi/q)
\right)\right\}}{n^{ps}}+\lim_{s\to 1/p}\,(2s-1)\zeta(ps)\,\pi\,{\rm
csc}\left(\frac{\pi}q \right)\\&=&0+\pi\,{\rm csc}\left(\frac{\pi}q
\right).\end{eqnarray*} Setting $r:=(p-1)s,s \to 1/p^+$ we check
that $\zeta(ps)^{1/p}\,\zeta(qr)^{1/q}=\zeta(ps)$ and hence the best
constant in (\ref{wittenp}) is the one given.

To recapitulate in terms of the celebrated infinite \emph{Hilbert
matrix,} $\mathcal{H}_0:=\left\{1/(m+n)\right\}_{m,n=1}^\infty$,
\cite[pp. 250--252]{exp2}, we have actually proven:

\begin{Thm} \label{hmat} Let $1<p,q < \infty$ be given with $1/p+1/q=1$.
The  Hilbert matrice $\mathcal{H}_0$  determines a bounded linear
mappings from the sequence space $\ell^p$ to itself such that
$$\|\mathcal{H}_0\|_{p,p}=\lim_{s \to
1/p}\frac{\mathcal{W}(s,(p-1)s,1)}{\zeta(ps)}=\pi\,{\rm
csc}\left(\frac{\pi}p \right).$$
\end{Thm}
\noindent{\bf Proof.}  Appealing to the isometry between
$(\ell^p)^*$ and $\ell^q$,  and given the evaluation $\mathcal{L}_p$
above, we directly compute the operator norm of $\mathcal{H}_0$ as
\begin{eqnarray*}\|\mathcal{H}_0\|_{p,p} = \sup_{\|x\|_p=1} \|\mathcal{H}_0
x\|_p=\sup_{\|y\|_q=1}\sup_{\|x\|_p=1} \langle \mathcal{H}_0 x,
y\rangle=\pi\,{\rm csc}\left(\frac{\pi}p \right).
\end{eqnarray*}
 \qed

A delightful  operator-theoretic introduction to the Hilbert matrix
$\mathcal{H}_0$ is given by Choi in his Chauvenet prize winning
article \cite{choi}.

One  may also study the corresponding behaviour of Hardy's
inequality (\ref{h-ineq}). For example, setting $a_n:=1/n$ in
(\ref{h-ineq})  and denoting $H_n:=\sum_{k=1}^n 1/k$ yields
$$\sum_{n=1}^\infty \left(\frac{H_n}{n}\right)^p \le
\left(\frac{p}{p-1}\right)^p\,\zeta(p).$$ Application of the
integral test and the evaluation $$ \int_1^\infty\,\left(\frac{\log
x}{x}\right)^p\,dx = \frac{\Gamma  \left( 1+p \right)}{ \left( p-1
\right) ^{p+1}},$$ for $p>1$ easily shows the constant is again best
possible.

\section{A Stirling
 Number Generating Function}\label{sect:Stirling}

Following~\cite{Butzer}, we begin with the integral
representation~\eqref{ZetaStirling} of \S\ref{sect:single2}. In
light of the expansion
\[
   \frac{(-1)^m}{m!}\log^m(1-x) = \sum_{n=0}^\infty u(n,m)
   \frac{x^n}{n!},\qquad 0\le m\in\Z,
\]
in terms of the unsigned Stirling numbers of the first kind (also
referred to as the Stirling cycle numbers in~\cite{GKP}), we have
\[
   \zeta(m+1) = \int_0^1 \bigg\{\sum_{n=1}^\infty u(n,m)\frac{x^n}{n!}\bigg\}\frac{dx}{x}
    = \sum_{n=1}^\infty \frac{u(n,m)}{n!\, n},\qquad 1\le m\in\Z.
\]
Telescoping the known recurrence
\begin{equation}\label{StirlingRecur}
   u(n,m) = u(n-1,m-1)+(n-1)u(n-1,m),\qquad 1\le m\le n,
\end{equation}
yields
\begin{equation}\label{StirlingSumRecur}
   u(n,m) = (n-1)! \left\{ \delta_{m,1}+ \sum_{j=1}^{n-1}
   \frac{u(j,m-1)}{j!}\right\}.
\end{equation}
Iterating this gives the representation
\[
   \zeta(m+1) = \zeta(2,\{1\}^{m-1}), \qquad 1\le m\in\Z,
\]
the $m=2$ case of which is~\eqref{z21}.  See also $n=0$
in~\eqref{DrinDuality} below.   \qed

For the alternating case, we begin by writing the
recurrence~\eqref{StirlingRecur} in the form
\[
   u(n+1,k) +(j-n)u(n,k) = u(n,k-1) + j\, u(n,k).
\]
Following~\cite{Butzer}, multiply both sides by
$(-1)^{n+k+1}j^{k-m-1}/(j-n)_n$, where $1\le n\le j-1$ and $k,
m\in\Z^{+}$, yielding
\begin{multline*}
   (-1)^k \left\{\frac{(-1)^{n+1}\,u(n+1,k)}{(j-n)_n} - \frac{(-1)^n
   \,u(n,k)}{(j-n+1)_{n-1}}\right\} j^{k-m-1}\\
   = \frac{(-1)^n}{(j-n)_n}\left\{(-1)^{k-1}\,u(n,k-1)j^{k-m-1}
   - (-1)^k \,u(n,k) j^{k-m}\right\}.
\end{multline*}
Now sum on $1\le k\le m$ and $1\le n\le j-1$, obtaining
\[
   \sum_{k=1}^m \frac{(-1)^{k+j} \,u(j,k)}{j!\,j^{m-k}} -\frac1{j^m}
   = \frac{(-1)^{m+1}}{(j-1)!}\sum_{n=m}^{j-1}(-1)^n
   (j-n-1)!\,u(n,m).
\]
Finally, sum on $j\in\Z^{+}$ to obtain
\[
   \zeta(m) = \sum_{k=1}^m \sum_{j=k}^\infty
   \frac{(-1)^{k+j}\,u(j,k)}{j!\,j^{m-k}}
   + \sum_{n=m}^\infty (-1)^{n+m}\,u(n,m)\sum_{j=n+1}^\infty
   \frac{(j-1-n)!}{(j-1)!}.
\]
Noting that
\[
  \sum_{j=n+1}^\infty \frac{(j-1-n)!}{(j-1)!}
  = \sum_{k=0}^\infty \frac{k!}{(k+n)!}
  = \frac1{n!} \;{}_2F_1(1,1;n+1;1)
  = \frac{1}{(n-1)!\,(n-1)},
\]
we find that
\[
   \zeta(m) = \sum_{k=1}^m \sum_{j=k}^\infty
   \frac{(-1)^{j+k}\,u(j,k)}{j!\,j^{m-k}} + \sum_{n=m}^\infty
   \frac{(-1)^{n+m}\,u(n,m)}{(n-1)!\,(n-1)}.
\]
Now employ the recurrence~\eqref{StirlingRecur} again to get
\begin{align}
   \zeta(m) &= \sum_{k=1}^{m-2}\sum_{j=k}^\infty
   \frac{(-1)^{j+k}\,u(j,k)}{j!\,j^{m-k}}+\sum_{j=m-1}^\infty
   \frac{(-1)^{j+m-1}\,u(j,m-1)}{j!\,j}+\sum_{j=m}^\infty
   \frac{(-1)^{j+m}\,u(j,m)}{j!}\nonumber\\
   &\qquad + \sum_{n=m}^\infty \frac{(-1)^{n+m}\,u(n-1,m)}{(n-1)!}
   + \sum_{n=m}^\infty
   \frac{(-1)^{n+m}\,u(n-1,m-1)}{(n-1)!\,(n-1)}\nonumber\\
   &=\sum_{k=1}^{m-2}\sum_{j=k}^\infty
   \frac{(-1)^{j+k}\,u(j,k)}{j!\,j^{m-k}}+2\sum_{j=m-1}^\infty
   \frac{(-1)^{j+m-1}\,u(j,m-1)}{j!\,j}.
\label{Butzer2bar}
\end{align}
Using~\eqref{StirlingSumRecur} again, we find that the case $m=3$
of~\eqref{Butzer2bar} gives
\begin{align*}
   \zeta(3) &= \sum_{j=1}^\infty
   \frac{(-1)^j\,u(j,1)}{j!\,j^2}+2\sum_{j=2}^\infty \frac{(-1)^j
   \,u(j,2)}{j!\,j}\\
   &= \sum_{j=1}^\infty \frac{(-1)^{j+1}}{j^3}+2\sum_{j=2}^\infty
  \frac{(-1)^j}{j!\,j}(j-1)!\sum_{k=1}^{j-1}\frac{u(k,1)}{k!}\\
   &= \sum_{j=1}^\infty \frac{(-1)^{j+1}}{j^3}+2\sum_{j=2}^\infty
   \frac{(-1)^j}{j^2}\sum_{k=1}^{j-1}\frac{1}{k}\\
   &= 2\zeta(\overline2,1)-\zeta(\overline3),
\end{align*}
which easily rearranges to give \eqref{dejavu}, shown in
\S\ref{sect:telescope} to be trivially equivalent to~\eqref{2bar1}.
\qed

\section{Polylogarithm Identities}\label{sect:dilog}

\subsection{Dilogarithm and Trilogarithm}\label{sect:dilog3} Consider the power series
\[
   J(x):= \zeta_x(2,1)=\sum_{n>k>0} \frac{x^n}{n^2 k},
   \qquad 0\le x\le 1.
\]
In light of~\eqref{iterint}, we have
\[
   J(x) = \int_0^x \frac{dt}{t}\int_0^t
   \frac{du}{1-u}\int_0^v\frac{dv}{1-v}
   = \int_0^x\frac{\log^2(1-t)}{2t}\,dt.
\]
The computer algebra package {\sc Maple} readily evaluates
\begin{equation}\label{maple}
   \int_0^x\frac{\log^2(1-t)}{2t}\,dt
   =\zeta(3)+\frac{1}{2} \log^2(1-x) \log(x)
       + \log (1-x){\rm Li}_2(1-x)-{\rm Li}_3(1-x)
\end{equation}
where
\[
   {\rm Li}_s(x) := \sum_{n=1}^\infty \frac{x^n}{n^s}
\]
is the classical polylogarithm~\cite{Lewin1,Lewin2}.  (One can
also readily verify the identity~\eqref{maple}  by differentiating
both sides by hand, and then checking~\eqref{maple} trivially
holds as $x\to0+$.  See also~\cite[p.\ 251, Entry 9]{Berndt1}.)
Thus,
\[
  J(x) =\zeta(3)+\frac{1}{2} \log^2(1-x) \log(x)
       + \log (1-x){\rm Li}_2(1-x)-{\rm Li}_3(1-x).
\]
Letting $x\to1-$ gives~\eqref{z21} again. \qed

In~\cite[p.\ 251, Entry 9]{Berndt1}, we also find that
\begin{equation}
   J(-z)+J(-1/z) = -\tfrac16\log^3 z-\mathrm{Li}_2(-z)\log
   z+\mathrm{Li}_3(-z)+\zeta(3)\label{Jinversion}
\end{equation}
and
\begin{equation}
   J(1-z) = \tfrac12\log^2 z\log(z-1)-\tfrac13\log^3
   z-\mathrm{Li}_2(1/z)\log
   z-\mathrm{Li}_3(1/z)+\zeta(3)\label{Jreflection}.
\end{equation}
Putting $z=1$ in~\eqref{Jinversion} and employing the well-known
dilogarithm evaluation~\cite[p.\ 4]{Lewin1}
\[
   \mathrm{Li}_{2}(-1) = \sum_{n=1}^\infty \frac{(-1)^n}{n^2} =
   -\frac{\pi^2}{12}
\]
gives~\eqref{2bar1}.  Putting $z=2$ in~\eqref{Jreflection} and
employing the dilogarithm evaluation~\cite[p.\ 6]{Lewin1}
\[
   \mathrm{Li}_{2}\bigg(\frac12\bigg) = \sum_{n=1}^\infty \frac1{n^2\, 2^n} =
   \frac{\pi^2}{12}-\frac12\log^2 2
\]
and the trilogarithm evaluation~\cite[p.\ 155]{Lewin1}
\[
   \mathrm{Li}_{3}\bigg(\frac12\bigg) = \sum_{n=1}^\infty \frac{1}{n^3\, 2^n} =
   \frac78\zeta(3)-\frac{\pi^2}{12}\log 2+\frac16\log^3 2
\]
gives~\eqref{2bar1} again. \qed

Finally, as in~\cite[Lemma~10.1]{BBBLa}, differentiation shows
that
\begin{equation}
\label{feq-J}
   J(-x) = -J(x)+\frac14J(x^2)+J\!\left(\frac{2x}{x+1}\right)
         -\frac18J\!\left(\frac{4x}{(x+1)^2}\right).
\end{equation}
Putting~\cite[Theorem~10.3]{BBBLa} $x=1$ gives $8J(-1)=J(1)$
immediately, i.e.~\eqref{2bar1}. \qed

In~\cite{BBBLa}, it is noted that once the component functions
in~\eqref{feq-J} are known, the coefficients can be deduced by
computing each term to high precision with a common transcendental
value of $x$ and then employing a linear relations finding
algorithm. We note here a somewhat more satisfactory method for
arriving at~\eqref{feq-J}.

First, as in \S\ref{sect:iterint1} one must determine the
fundamental transformations~\eqref{identity}--\eqref{quadLanden}.
While this is not especially difficult, as the calculations are
somewhat lengthy, we do not include them here.  By performing
these transformations on the function $J(x)$, one finds that
\begin{alignat*}{2}
   J(x) & = \int_0^x ab^2,\qquad\qquad
        & J\bigg(\frac{2x}{1+x}\bigg) &= \int_0^x (a+c)(b-c)^2,\\
   J(-x) & = \int_0^x ac^2, \qquad\qquad &{} &{}\\
   J(x^2) &= \int_0^x 2a(b+c)^2,  \qquad\qquad
      & J\bigg(\frac{4x}{(1+x)^2}\bigg)  &= \int_0^x (a+2c)4(b-c)^2.
\end{alignat*}
It now stands to reason that we should seek rational numbers
$r_1$, $r_2$, $r_3$ and $r_4$ such that
\[
   ac^2 = r_1 ab^2 + 2r_2 \, a(b+c)^2 + r_3(a+c)(b-c)^2
   +r_4(a+2c)4(b-c)^2
\]
is an identity in the non-commutative polynomial ring $\Q\langle
a,b,c\rangle$.  The problem of finding such rational numbers reduces
to solving a finite set of linear equations.  For example, comparing
coefficients of the monomial $ab^2$ tells us that
$r_1+2r_2+r_3+4r_4=0$.  Coefficients of other monomials give us
additional equations, and we readily find that $r_1=-1$, $r_2=1/4$,
$r_3=1$ and $r_4=-1/8$, thus proving~\eqref{feq-J} as expected.

\subsection{Convolution of Polylogarithms}\label{sect:dilogp}
Motivated by~\cite{Boy2001,Boy2002}, for real $0<x<1$ and integers
$s$ and $t$, consider
\begin{align*}
   T_{s,t}(x) &:= \sum_{\substack{m,n=1\\ m\ne n}}^\infty
   \frac{x^{n+m}}{n^s\, m^t(m-n)}
   =\sum_{\substack{m,n=1\\ m\ne n}}^\infty
   \frac{x^{n+m}(m-n+n)}{n^s\, m^{t+1}(m-n)}\\
   &=\sum_{\substack{m,n=1\\ m\ne n}}^\infty
   \frac{x^{n+m}}{n^s \, m^{t+1}}+\sum_{\substack{m,n=1\\ m\ne n}}^\infty
   \frac{x^{n+m}}{n^{s-1}\, m^{t+1}(m-n)}\\
   &=\sum_{n=1}^\infty \frac{x^n}{n^s}\sum_{m=1}^\infty
   \bigg(\frac{x^m}{m^{t+1}}-\frac{x^n}{n^{t+1}}\bigg)+T_{s-1,t+1}(x)\\
   &=\mathrm{Li}_s(x)\mathrm{Li}_{t+1}(x)-\mathrm{Li}_{s+t+1}(x^2)
   +T_{s-1,t+1}(x).
\end{align*}
Telescoping this gives
\begin{align}
   T_{s,t}(x) &= T_{0,s+t}(x)-s\,\mathrm{Li}_{s+t+1}(x^2)
   +\sum_{j=1}^s
   \mathrm{Li}_j(x)\mathrm{Li}_{s+t+1-j}(x),
   \qquad 0\le s\in\Z.\nonumber\\
\intertext{With $t=0$, this becomes}
   T_{s,0}(x) &= T_{0,s}(x) - s\,\mathrm{Li}_{s+1}(x^2)
              +\sum_{j=1}^{s}\mathrm{Li}_j(x)\mathrm{Li}_{s+1-j}(x),
              \qquad 0\le s\in\Z.\nonumber\\
\intertext{But for any integers $s$ and $t$, there holds}
   T_{s,t}(x) &= \sum_{\substack{m,n=1\\ m\ne n}}^\infty
   \frac{x^{n+m}}{n^t m^s(m-n)}
   = -\sum_{\substack{m,n=1\\ m\ne n}}^\infty
   \frac{x^{n+m}}{m^s n^t(n-m)}
   = - T_{s,t}(x).\nonumber\\
\intertext{Therefore,}
   T_{s,0}(x) &= \frac12 \sum_{j=1}^s \mathrm{Li}_j(x)\mathrm{Li}_{s+1-j}(x)
              - \frac{s}{2}\,\mathrm{Li}_{s+1}(x^2),
              \qquad 0\le s\in\Z.\label{Tp0}\\
\intertext{On the other hand,}
   T_{s,0}(x) &= \sum_{n=1}^\infty\frac{x^n}{n^s}
   \sum_{\substack{m=1\\m\ne n}}^\infty\frac{x^m}{m-n}
   = \sum_{n=1}^\infty\frac{x^{2n}}{n^s}\sum_{m=n+1}^\infty
   \frac{x^{m-n}}{m-n} - \sum_{n=1}^\infty
   \frac{x^n}{n^s}\sum_{m=1}^{n-1}\frac{x^m}{n-m}\nonumber\\
   &= \mathrm{Li}_s(x^2)\mathrm{Li}_1(x)-\sum_{n=1}^\infty
   \frac{x^n}{n^s}\sum_{j=1}^{n-1}\frac{x^{n-j}}{j}.\nonumber\\
\intertext{Comparing this with~\eqref{Tp0} gives}
   \sum_{n=1}^\infty\frac{x^n}{n^s}\sum_{j=1}^{n-1}\frac{x^{n-j}}{j}
    & =\frac{s}2\,\mathrm{Li}_{s+1}(x^2)
      -\left[\mathrm{Li}_s(x)-\mathrm{Li}_s(x^2)\right]\mathrm{Li}_1(x)
      -\frac12\sum_{j=2}^{s-1}\mathrm{Li}_j(x)\mathrm{Li}_{s+1-j}(x),\label{TakeLimit}
\end{align}
where in~\eqref{TakeLimit} and what follows, we now require $2\le
s\in\Z$ because the terms $j=1$ and $j=s$ in the sum~\eqref{Tp0}
were separated, and assumed to be distinct.

Next, note that if $n$ is a positive integer and $0<x<1$, then
\[
   1-x^n = (1-x)\sum_{j=0}^{n-1} x^j < (1-x)n.
\]
Thus, if $2\le s\in\Z$ and $0<x<1$, then
\begin{align*}
   0<\left[\mathrm{Li}_s(x)-\mathrm{Li}_s(x^2)\right]\mathrm{Li}_1(x)
   &= \mathrm{Li}_1(x)\sum_{n=1}^\infty \frac{x^n(1-x^n)}{n^s}
   < (1-x)\mathrm{Li}_1(x)\sum_{n=1}^\infty \frac{x^n}{n^{s-1}}\\
   &< (1-x)\log^2(1-x).
\end{align*}
Since the latter expression tends to zero in the limit as $x\to
1-$, taking the limit in~\eqref{TakeLimit} gives
\[
   \zeta(s,1) = \frac12\, s\,\zeta(s+1)-\frac12\sum_{j=1}^{s-2}
   \zeta(j+1)\zeta(s-j),
   \qquad 2\le s\in\Z,
\]
which is~\eqref{EulerReduction}. \qed

\section{Fourier Series}\label{sect:fourier}

The Fourier expansions
\[
   \sum_{n=1}^{\infty}\frac{\dst
   \sin(nt)}{n}=\frac{\dst\pi-t}{2}\qquad\text{and}\qquad
   \sum_{n=1}^{\infty} \frac{\dst\cos(nt)}{n}= -\log|2\sin(t/2)|
\]
are both valid in the open interval $0<t<2\pi$.   Multiplying
these together, simplifying, and doing a partial fraction
decomposition gives
\begin{align}
   \sum_{n=1}^{\infty}\frac{\sin(nt)}{n} \sum_{k=1}^{n-1} \frac 1 k
   & = \frac12\sum_{n=1}^\infty\frac{\sin(nt)}{n}
       \sum_{k=1}^{n-1}\bigg(\frac1k+\frac1{n-k}\bigg)
   = \frac12\sum_{n>k>0} \frac{\sin(nt)}{k(n-k)}\nonumber\\
   &= \frac12\sum_{m,n=1}^\infty \frac{\sin(m+n)t}{mn}= \sum_{m,n=1}^\infty \frac{\sin(mt)\cos(nt)}{mn}\nonumber\\
   &= -\frac{\pi-t}2 \log \left |2\sin(t/2) \right|, \label{sin-ser}
\intertext{again for $0<t<2\pi$.  Integrating (\ref{sin-ser}) term
by term yields}
   \sum_{n=1}^{\infty}\frac{\cos(n \theta)}{n^2} \sum_{k=1}^{n-1}\frac1k
   &= \zeta(2,1)+\frac12\int_0^\theta (\pi -t)\log \left
   |2\sin(t/2)\right|\,dt, \label{z21f}\\
\intertext{valid for $0\le \theta\le 2\pi$.  Likewise for $0\le
\theta\le 2\pi$,}
   \sum_{n=1}^{\infty}\frac{ \cos(n \theta)}{n^3}
   &= \zeta(3)+\int_0^\theta (\theta-t)\log \left |2\sin(t/2)
      \right|\,dt.\label{z3f}
\end{align}
Setting $\theta=\pi$ in~\eqref{z21f}~and~\eqref{z3f} produces
\[
   \zeta(2,1)-\zeta(\overline{2},1) = -\frac{1}{2}\,\int_0^{\pi}(\pi-t)
   \log \left |2\sin(t/2) \right| \,dt =
   \frac{\zeta(3)-\zeta(\overline{3})}2.
\]
In light of~\eqref{z21}, this implies
\[
   \zeta(\overline2,1)=\frac{\zeta(3)+\zeta(\overline3)}{2}=
   \frac12\sum_{n=1}^\infty\frac{1+(-1)^n}{n^3}=\sum_{m=1}^\infty\frac{1}{(2m)^3}
   =\frac18\zeta(3),
\]
which is~\eqref{2bar1}. \qed

Applying Parseval's equation to~\eqref{sin-ser} gives
(via~\cite{BB,BBG,Flaj}) the integral evaluation
\[
   \frac{1}{4\pi} \int_{0}^{2\pi} (\pi-t)^2 \log^2(2 \sin(t/2))\, dt =
   \sum_{n=1}^{\infty} \frac{H_n^2}{(n+1)^2} = \frac{11}4\, \zeta(4).
\]
A reason for valuing such integral representations is that they are
frequently easier to use numerically.

\section{Further Generating Functions}\label{sect:gfs}

\subsection{Hypergeometric Functions}\label{sect:hyper}

Note that in the notation of~\eqref{mzvdef}, $\zeta(2,1)$ is the
coefficient of $xy^2$ in
\begin{equation}\label{xyy}
   G(x,y) :=
   \sum_{m=0}^\infty\sum_{n=0}^\infty x^{m+1}y^{n+1}\zeta(m+2,\us^n)
   = y\sum_{m=0}^\infty x^{m+1}\sum_{k=1}^\infty
      \frac1{k^{m+2}}\prod_{j=1}^{k-1} \left(1+\frac{y}j\right).
\end{equation}
Now recall the notation $(y)_k := y(y+1)\cdots(y+k-1)$ for the
rising factorial with $k$ factors.  Thus,
\[
   \frac{y}{k}\prod_{j=1}^{k-1} \left(1+\frac{y}j\right)
   = \frac{(y)_k}{k!}.
\]
Substituting this into~\eqref{xyy}, interchanging order of
summation, and summing the resulting geometric series yields the
hypergeometric series
\begin{align*}
   G(x,y)
   &= \sum_{k=1}^\infty\frac{(y)_k}{k! }\left(\frac{x}{k-x}\right)
   = -\sum_{k=1}^\infty\frac{(y)_k(-x)_k}{k! (1-x)_k}
   = 1-{}_2F_1\left(\begin{array}{cc} y, -x\\ 1-x\end{array}\bigg|1\right).
\end{align*}
 But, Gauss's summation theorem for the hypergeometric
function~\cite[p.\ 557]{AS}~\cite[p.\ 2]{Bailey} and the power
series expansion for the logarithmic derivative of the gamma
function \cite[p.\ 259]{AS} imply that
\[
   {}_2F_1\left(\begin{array}{cc} y, -x\\ 1-x\end{array}\bigg|1\right)
   =\frac{\G(1-x)\G(1-y)}{\G(1-x-y)}
   =\exp\bigg\{\sum_{k=2}^\infty\left(x^k+y^k-(x+y)^k\right)
   \frac{\zeta(k)}{k}\bigg\}.
\]
Thus, we have derived the generating function equality~\cite{BBB}
(see~\cite{DBq} for a $q$-analog)
\begin{equation}\label{DrinGF}
   \sum_{m=0}^\infty\sum_{n=0}^\infty x^{m+1}y^{n+1}\zeta(m+2,\us^n)
   = 1-\exp\bigg\{\sum_{k=2}^\infty
     \left(x^k+y^k-(x+y)^k\right)\frac{\zeta(k)}{k}\bigg\}.
\end{equation}
Extracting coefficients of $xy^2$ from both sides
of~\eqref{DrinGF} yields~\eqref{z21}. \qed

The generalization~\eqref{EulerReduction} can be similarly
derived: extract the coefficient of $x^{m-1}y^2$ from both sides
of~\eqref{DrinGF}. In fact, it is easy to see that~\eqref{DrinGF}
provides a formula for $\zeta(m+2,\us^n)$ for all nonnegative
integers $m$ and $n$ in terms of sums of products of values of the
Riemann zeta function at the positive integers. In particular,
Markett's formula~\cite{Mark} (cf.\ also~\cite{BBG}) for
$\zeta(m,1,1)$ for positive integers $m>1$ is most easily obtained
in this way. Noting symmetry between $x$ and $y$ in~\eqref{DrinGF}
gives Drinfeld's duality formula~\cite{Drin}
\begin{equation}\label{DrinDuality}
   \zeta(m+2,\us^n) = \zeta(n+2,\us^m)
\end{equation}
for non-negative integers $m$ and $n$, a special case of the more
general duality formula~\eqref{duality}.  Note that~\eqref{z21} is
just the case $m=n=0$.

Similarly~\cite[2.1b]{Chu} equating coefficients of $xy^2$ in
Kummer's summation theorem~\cite[p.\ 53]{Kummer}~\cite[p.\
9]{Bailey}
\[
   {}_2F_1\left(\begin{array}{cc}x,y\\1+x-y\end{array}\bigg|-1\right)
   = \frac{\Gamma(1+x/2)\Gamma(1+x-y)}{\Gamma(1+x)\Gamma(1+x/2-y)}
\]
yields~\eqref{2bar1}.

\subsection{A Generating Function for Sums}\label{sect:SumGF}

The identity~\eqref{z21} can also be recovered by setting $x=0$ in
the following result:

\begin{Thm}[Theorem 1 of~\cite{DJD}]\label{thm:SumGF}
If $x$ is any complex number not equal to a
positive integer, then
\[
   \sum_{n=1}^\infty\frac1{n(n-x)}\sum_{m=1}^{n-1}\frac1{m-x}
   =\sum_{n=1}^\infty\frac1{n^2(n-x)}.
\]
\end{Thm}

\noindent{\bf Proof.} Fix $x\in\C\setminus\Z^{+}$. Let $S$ denote
the left hand side.   By partial fractions,
\begin{align*}
   S & 
   =\sum_{n=1}^\infty \sum_{m=1}^{n-1}\bigg(\frac{1}{n(n-m)(m-x)}
   -\frac{1}{n(n-m)(n-x)}\bigg)\\
   &=\sum_{m=1}^\infty\frac1{m-x}\sum_{n=m+1}^\infty\frac1{n(n-m)}
   -\sum_{n=1}^\infty \frac1{n(n-x)}\sum_{m=1}^{n-1}\frac1{n-m}\\
  &=\sum_{m=1}^\infty\frac1{m(m-x)}\sum_{n=m+1}^\infty\bigg(\frac1{n-m}
   -\frac1n\bigg)-\sum_{n=1}^\infty\frac1{n(n-x)}\sum_{m=1}^{n-1}\frac1m.
\end{align*}
Now for fixed $m\in\Z^{+}$,
\begin{align*}
   \sum_{n=m+1}^\infty \bigg(\frac1{n-m}-\frac1n\bigg)
   &= \lim_{N\to\infty}\sum_{n=m+1}^N \bigg(\frac1{n-m}-\frac1n\bigg)
   = \sum_{n=1}^m\frac1n-\lim_{N\to\infty}\sum_{n=1}^m\frac1{N-n+1}\\
   &= \sum_{n=1}^m \frac1n,
\end{align*}
since $m$ is fixed. Therefore, we have
\begin{align*}
   S &=\sum_{m=1}^\infty\frac1{m(m-x)}\sum_{n=1}^m\frac1n
   - \sum_{n=1}^\infty\frac1{n(n-x)}\sum_{m=1}^{n-1}\frac1m
   =\sum_{n=1}^\infty\frac1{n(n-x)}\bigg(\sum_{m=1}^n\frac1m-
   \sum_{m=1}^{n-1}\frac1m\bigg)\\
   &=\sum_{n=1}^\infty\frac1{n^2(n-x)}.
\end{align*}
\eop

Theorem~\ref{thm:SumGF} is in fact equivalent to the sum
formula~\cite{Gran,Ohno}
\begin{equation}\label{sum}
  \sum_{\substack{\sum a_i = s\\a_i\ge 0}}
  \zeta(a_1+2,a_2+1,\dots,a_r+1) = \zeta(r+s+1),
\end{equation}
valid for all integers $s\ge 0$, $r\ge 1$, and which generalizes
Theorem~\ref{briggs}~\eqref{sumdepth2} to arbitrary depth.   The
identity~\eqref{z21} is simply the case $r=2$, $s=0$. A $q$-analog
of the sum formula~\eqref{sum} is derived as a special case of
more general results in~\cite{DBq}.  See also~\cite{DBqSum}.

\subsection{An Alternating Generating Function}\label{sect:SumAGF}
An alternating counterpart to Theorem~\ref{thm:SumGF} is given
below.

\begin{Thm}{\rm (Theorem 3 of \cite{DJD}).} For all non-integer $x$
\begin{align*}
   \sum_{n=1}^\infty \frac{(-1)^n}{n^2-x^2}
      \bigg\{H_n +\sum_{n=1}^\infty \frac{x^2}{n(n^2-x^2)}\bigg\}
   &= \sum_{n=1}^\infty \frac{(-1)^n}{n^2-x^2}\bigg\{\psi(n)-\psi(x)
      - \frac{\pi}{2}\cot(\pi x) -\frac1{2x}\bigg\}\\
   &= \sum _{o>0\,\mbox{{\small odd}}}^{\infty }{\frac {1}{
o \left( o^{2 }-x^2 \right) }}+\sum _{n=1}^{\infty }{\frac {
\left( -1 \right) ^{n}\,n}{ \left( {n}^{2}-{x}^{2} \right)
^{2}}}\\
   &=\sum _{e>0\,\mbox{{\small even}}}^{\infty }{\frac {e}{
\left( {x}^{2}-e^{2} \right) ^{2}}}-{x}^{2}\sum
_{o>0\,\mbox{{\small odd}}}^{\infty }{\frac {1}{ o
 \left( {x}^{2}- o^{2} \right) ^{2}}}.
 \end{align*}
 \end{Thm}

Setting $x=0$ reproduces \eqref{2bar1} in the form
$\zeta(\overline{2},1)=\sum _{n>0}^{\infty }(2n)^{-3}$. We record
that
\[
   \sum _{n=1}^{\infty}\frac{(-1)^n}{n^2-x^2}
   =\frac1{2 x^2}-{\frac {\pi }{2 x\sin(\pi x)}},
\]
while
\begin{multline*}
   \sum_{n=1}^\infty \frac{(-1)^n}{n^2-x^2}
        \bigg\{\psi(n)-\psi(x)-\frac\pi2\cot(\pi
        x)-\frac1{2x}\bigg\}
  =\sum_{n=1}^\infty \frac{(-1)^n}{n^2-x^2}\bigg\{H_n
    +\sum_{n=1}^\infty \frac{x^2}{n(n^2-x^2)} \bigg\}\\
   = \sum_{n=1}^\infty\frac{1}{(2n-1)((2n-1)^2-x^2)}
    +\sum_{n=1}^\infty \frac{n (-1)^n}{(n^2-x^2)^2}.
\end{multline*}

\subsection{The Digamma Function} \label{susbsect:DPsiGF}
Define an auxiliary function $\Lambda$ by
\begin{align*}
   x \Lambda(x):= \tfrac12\psi'(1-x)
       -\tfrac12\left(\psi(1-x)+\gamma\right)^2
       -\tfrac12\zeta(2).
\end{align*}
We note, but do not use, that
\[
   x \Lambda(x)=\frac12\int_0^\infty
   \frac{t \left(e^{-t}+e^{-t (1-x)}\right)}{1-e^{-t}}\,dt
   -\frac12\left(\int_0^\infty
      \frac{e^{-t}-e^{-t ( 1-x)}}{1-{e^{-t}}}\,dt\right)^2-\zeta(2).
\]
It is easy to verify that
\begin{align}\label{2ids}
   \psi(1-x)+\gamma &=
   \sum_{n=1}^\infty\frac{x}{n(x-n)},\nonumber\\
   \psi'(1-x)-\zeta(2)
   & =\sum_{n=1}^\infty\bigg(\frac{1}{(x-n)^2}-\frac1{n^2}\bigg)
   =\sum_{n=1}^\infty \frac{2nx-x^2}{n^2(n-x)^2},
\end{align}
and
\[
   \sum_{n=0}^\infty\zeta(n+2,1)x^n
   =\sum_{n=1}^\infty\frac1{n(n-x)}\sum_{m=1}^{n-1}\frac1m.
\]
Hence,
\[
   \Lambda(x)
    =\sum_{n=1}^\infty\frac{1}{n^2(n-x)}
    -x\sum_{n=1}^\infty\frac1{n(n-x)}\sum_{m=1}^{n-1}\frac{1}{m(m-x)}.
\]
Now,
\[
    \sum_{n=1}^\infty\frac{1}{n^2(n-x)} -
    x\sum_{n=1}^\infty\frac{1}{n(n-x)}\sum_{m=1}^{n-1}\frac{1}{m(m-x)}
    =\sum_{n=1}^\infty\frac1{n(n-x)}\sum_{m=1}^{n-1}\frac1m
\]
is directly equivalent to Theorem~\ref{thm:SumGF} of
\S\ref{sect:SumGF}---see \cite[Section 3]{DJD}---and we have
proven $$\Lambda(x)= \sum_{n=0}^\infty \zeta(n+2,1)\,x^n,$$ so
that comparing coefficients yields
  yet
another proof of Euler's reduction~\eqref{EulerReduction}.  In
particular, setting $x=0$ again produces~\eqref{z21}. \qed

\subsection{The Beta Function}

Recall that the beta function is defined for positive real $x$ and
$y$ by
\[
   B(x,y) := \int_0^1 t^{x-1}(1-t)^{y-1}\,dt
           = \frac{\Gamma(x)\Gamma(y)}{\Gamma(x+y)}.
\]
We begin with the following easily obtained generating function:
\[
   \sum_{n=1}^\infty t^n H_n  = -\frac{\log(1-t)}{1-t}.
\]
For $m\ge 2$ the Laplace integral~\eqref{Laplace} now gives
\begin{align}\label{zm1=b1}
   \zeta(m,1) & = \frac{(-1)^m}{(m-1)!}
      \int_0^1\frac{\log^{m-1}(t)\log(1-t)}{1-t}\, dt\nonumber\\
   & = \frac{(-1)^m}{2(m-1)!}
       \int_0^1 (m-1) \log^{m-2}(t)\log^2(1-t)\,\frac{dt}{t}\nonumber\\
   & = \frac{(-1)^m}{2(m-2)!}\, b_1^{(m-2)}(0),
\end{align}
where
\[
   b_1(x) := \left. \frac{\partial^{2}}{\partial y^2} B(x,y)
             \right|_{y=1}
           = 2\Lambda(-x)
\]
(cf.~\S\ref{susbsect:DPsiGF}). Since
\[
  \frac{\partial^2}{\partial y^2} B(x,y)
  = B(x,y) \left[(\psi(y)-\psi(x+y))^2 + \psi'(y) - \psi'(x+y)\right],
\]
we derive
\[
   b_1(x) = \frac{(\psi(1) - \psi(x+1))^2 + \psi'(1) - \psi'(x+1)}{x}.
\]
Now observe that from~\eqref{zm1=b1},
\begin{align*}
   \zeta(2,1) &= \frac12 \, b_1(0)
   = \lim_{x \downarrow0} \frac{(-\gamma - \psi(x+1))^2 }{2x}
    -\lim_{x \downarrow0} \frac{\psi'(x+1)-\psi'(1)}{2x}
   =-\frac12\,\psi''(1)\\
 &=\zeta(3).
\end{align*}
\qed

Continuing, from the following two identities, cognate to
\eqref{2ids},
\begin{align*}
   (-\gamma-\psi(x+1))^2
   & =  \bigg(\sum_{m=1}^\infty (-1)^m\zeta(m+1)\, x^m \bigg)^2\\
   & =  \sum_{m=1}^\infty (-1)^m
        \sum_{k=1}^{m-1}\zeta(k+1)\zeta(m-k+1)\, x^m,\\
   \zeta(2) -\psi'(x+1)
   & = \sum_{m=1}^\infty (-1)^{m+1} (m+1) \zeta(m+2)\, x^m,
\end{align*}
we get
\[
   2\sum_{m=2}^\infty (-1)^m \zeta(m,1) \,x^{m-2}
   =  \sum_{m=2}^\infty \frac{b_1^{(m-2)}(0)}{(m-2)!}\, x^{m-2}\,
   =  b_1(x)
\]
\[
   = \sum_{m=1}^\infty (-1)^{m-1} \left( (m+1) \zeta(m+2) -
\sum_{k=1}^{m-1} \zeta(k+1) \zeta(m-k+1)\right) x^{m-1},
\]
from which Euler's reduction~\eqref{EulerReduction}
follows---indeed this is close to Euler's original path.

Observe that~\eqref{zm1=b1} is especially suited to symbolic
computation. We also note  the pleasing identity
\begin{eqnarray}\psi'(x)=\frac{\Gamma''(x)}{\Gamma(x)}-\psi^2(x)\label{DPsiGF}
.\end{eqnarray} In some informal sense (\ref{DPsiGF})
generates~\eqref{EulerReduction}, but we have been unable to make
this  sense precise.

\section{A Decomposition Formula of Euler}\label{sect:parfracs}
For positive integers $s$ and $t$ and distinct non-zero real
numbers $\alpha$ and $x$, the partial fraction expansion
\begin{equation}\label{parfrac}
   \frac{1}{x^s(x-\alpha)^t}
   =
   (-1)^t\sum_{r=0}^{s-1}\binom{t+r-1}{t-1}\frac{1}{x^{s-r}\alpha^{t+r}}
   +\sum_{r=0}^{t-1}\binom{s+r-1}{s-1}
   \frac{(-1)^r}{\alpha^{s+r}(x-\alpha)^{t-r}}
\end{equation}
implies~\cite[p.\ 48]{Niels3}~\cite{Mark} Euler's decomposition
formula
\begin{multline}\label{niels}
   \zeta(s,t) =
   (-1)^t\sum_{r=0}^{s-2}\binom{t+r-1}{t-1}\zeta(s-r,t+r)
   +\sum_{r=0}^{t-2}(-1)^r\binom{s+r-1}{s-1}\zeta(t-r)\zeta(s+r)\\
   -(-1)^t\binom{s+t-2}{s-1}\big\{\zeta(s+t)+\zeta(s+t-1,1)\big\}.
\end{multline}
The depth-2 sum formula~\eqref{sumdepth2} is obtained by setting
$t=1$ in~\eqref{niels}. If we also set $s=2$, the
identity~\eqref{z21} results.  To derive~\eqref{niels}
from~\eqref{parfrac} we follow~\cite{Mark}, separating the last
term of each sum on the right hand side of~\eqref{parfrac},
obtaining
\begin{align*}
   \frac{1}{x^s(x-\alpha)^t} &=
   (-1)^t\sum_{r=0}^{s-2}\binom{t+r-1}{t-1}\frac{1}{x^{s-r}\alpha^{t+r}}
   +\sum_{r=0}^{t-2}\binom{s+r-1}{s-1}
   \frac{(-1)^r}{\alpha^{s+r}(x-\alpha)^{t-r}}\\
   &-(-1)^{t}\binom{s+t-2}{s-1}
   \frac{1}{\alpha^{s+t-1}}\bigg(\frac{1}{x-\alpha}
   -\frac1x\bigg).
\end{align*}
Now sum over all integers $0<\alpha <x<\infty$. \qed

Nielsen states~\eqref{parfrac} without proof~\cite[p.\ 48, eq.\
(9)]{Niels3}. Markett proves~\eqref{parfrac} by
induction~\cite[Lemma 3.1]{Mark}, which is the proof technique
suggested for the $\alpha=1$ case of~\eqref{parfrac}
in~\cite[Lemma 1] {BBG}. However, it is easy to
prove~\eqref{parfrac} directly by expanding the left hand side
into partial fractions with the aid of the residue calculus.
Alternatively, as in~\cite{DBqDecomp} note that~\eqref{parfrac} is
an immediate consequence of applying the partial derivative
operator
\[
   \frac{1}{(s-1)!}\bigg(-\frac{\partial}{\partial x}\bigg)^{s-1}
   \frac{1}{(t-1)!}\bigg(-\frac{\partial}{\partial y}\bigg)^{t-1}
\]
to the identity
\[
   \frac{1}{xy} = \frac{1}{(x+y)x} + \frac{1}{(x+y)y},
\]
and then setting $y=\alpha-x$.  This latter observation is
extended in~\cite{DBqDecomp} to establish a $q$-analog of another
of Euler's decomposition formulas for $\zeta(s,t)$.

\section{Equating Shuffles and
Stuffles}\label{sect:shtuff}

We begin with an informal argument. By the \emph{stuffle
multiplication} rule~\cite{BBBLa,BowBradSurvey,DBPrtn,DBq}
\begin{equation}\label{divergentstuff}
   \zeta(2)\zeta(1)= \zeta(2,1)+\zeta(1,2)+\zeta(3).
\end{equation}
On the other hand, the shuffle multiplication
rule~\cite{BBBLa,BBBLc,BowBradSurvey,BowBrad3,BowBradRyoo} gives
$ab\shuff b = 2abb+bab$, whence
\begin{equation}\label{divergentshuff}
   \zeta(2)\zeta(1) = 2\zeta(2,1)+\zeta(1,2).
\end{equation}
The identity~\eqref{z21} now follows immediately on
subtracting~\eqref{divergentstuff} from~\eqref{divergentshuff}.
\qed

Of course, this argument needs justification, because it involves
cancelling divergent series. To make the argument rigorous, we
introduce the multiple
polylogarithm~\cite{BBBLa,BowBrad1,BowBradSurvey}. For real $0\le
x\le 1$ and positive integers $s_1,\dots,s_k$ with $x=s_1=1$
excluded for convergence, define
\begin{equation}
   \zeta_x(s_1,\dots,s_k)
   := \sum_{n_1>\cdots>n_k>0}\; x^{n_1}\prod_{j=1}^k n_j^{-s_j}
   =\int \prod_{j=1}^k
   \bigg(\prod_{r=1}^{s_j-1} \frac{dt_r^{(j)}}{t_r^{(j)}}\bigg)
   \frac{dt_{s_j}^{(j)}}{1-t_{s_j}^{(j)}},
\label{iterint}
\end{equation}
where the integral is over the simplex
\[
x>t_1^{(1)}>\cdots>t_{s_1}^{(1)}>\cdots>t_1^{(k)}>\cdots>t_{s_k}^{(k
)}>0,
\]
and is abbreviated by
\begin{equation}
   \int_0^x \prod_{j=1}^k a^{s_j-1}b,
   \qquad a=\frac{dt}{t},\qquad b=\frac{dt}{1-t}.
\label{shortiterint}
\end{equation}
Then
\[
   \zeta(2)\zeta_x(1) =
   \sum_{n>0}\frac1{n^2}\sum_{k>0}\frac{x^k}{k}
   =
   \sum_{n>k>0}\frac{x^k}{n^2k}+\sum_{k>n>0}\frac{x^k}{kn^2}
   +\sum_{k>0}\frac{x^k}{k^3},
\]
and
\[
   \zeta_x(2)\zeta_x(1) = \int_0^x ab \int_0^x b = \int_0^x
   \left(2abb+bab\right) = 2\zeta_x(2,1)+\zeta_x(1,2).
\]
Subtracting the two equations gives
\[
   \big[\zeta(2)-\zeta_x(2)\big]\zeta_x(1)
   = \zeta_x(3)-\zeta_x(2,1)+\sum_{n>k>0}\frac{x^k-x^n}{n^2k}.
\]
We now take the limit as $x\to 1-.$ Uniform convergence implies
the right hand side tends to $\zeta(3)-\zeta(2,1)$.  That the left
hand side tends to zero follows immediately from the inequalities
\begin{align*}
   0\le x \big[\zeta(2)-\zeta_x(2)\big]\zeta_x(1)
   &= x\int_x^1\log(1-t)\log(1-x)\frac{dt}{t}\\
   & \le \int_x^1 \log^2(1-t)\,dt\\
   & =(1-x)\left\{1+(1-\log(1-x))^2\right\}.
\end{align*} \qed

The alternating case~\eqref{2bar1} is actually easier using this
approach, since the role of the divergent sum $\zeta(1)$ is taken
over by the conditionally convergent sum $\zeta(\overline 1)=-\log
2$.  By the stuffle multiplication rule,
\begin{align}
   \zeta(\overline 2)\zeta(\overline 1) &= \zeta(\overline 2,
   \overline 1)+\zeta(\overline 1, \overline 2) + \zeta(3)
   \label{convergentstuff1},\\
   \zeta(2)\zeta(\overline 1) &= \zeta(2, \overline
   1)+\zeta(\overline 1, 2) + \zeta(\overline
   3)\label{convergentstuff2}.
\end{align}
On the other hand, the shuffle multiplication rule gives $ac\shuff
c=2ac^2+cac$ and $ab\shuff c=abc+acb+cab$, whence
\begin{align}
   \zeta(\overline 2)\zeta(\overline 1) &= 2\zeta(\overline 2, 1)
   + \zeta(\overline 1, 2)\label{convergentshuff1},\\
   \zeta(2)\zeta(\overline 1) &= \zeta(2, \overline 1) +
   \zeta(\overline 2, \overline 1) + \zeta(\overline 1, \overline
   2)\label{convergentshuff2}.
\end{align}
Comparing~\eqref{convergentstuff1} with~\eqref{convergentshuff1}
and~\eqref{convergentstuff2} with~\eqref{convergentshuff2} yields
the two equations
\begin{align*}
   \zeta(\overline 2, \overline 1)
   &= \zeta(\overline 1, 2) +2\zeta(\overline 2, 1)
   -\zeta(\overline 1, \overline 2)
   -\zeta(3),\\
   \zeta(\overline 2, \overline 1) &= \zeta(\overline 1, 2)-\zeta(
   \overline 1, \overline 2) +\zeta(\overline 3).
\end{align*}
Subtracting the latter two equations yields $2\zeta(\overline 2,
1) = \zeta(3)+\zeta(\overline 3)$, i.e.~\eqref{dejavu}, which was
shown to be trivially equivalent to~\eqref{2bar1} in
\S\ref{sect:telescope}. \qed

\section{Conclusion}

There are doubtless other roads to Rome, and as indicated in the
introduction we should like to learn of them. We finish with the
three open questions we are most desirous of answers to.

\begin{itemize}
\item A truly combinatorial proof, perhaps of the form considered in
      \cite{BBBLc}.
\item A direct proof that the appropriate line integrals in
      sections~\ref{sect:Perron} and~\ref{sect:Dirichlet} evaluate
      to the appropriate multiple of $\zeta(3)$.
\item A proof of~\eqref{2bar1^n}, or at least some additional
cases of it.
\end{itemize}


\end{document}